\RequirePackage{fix-cm}
\documentclass[twocolumn]{svjour3}          
\smartqed  
\usepackage{graphicx}
\usepackage{nicefrac}
 \usepackage{mathptmx}      
\usepackage{hyperref}
\usepackage{color}
\usepackage{amsmath}
\usepackage{amssymb}
\usepackage{eucal}
\usepackage{booktabs}
\usepackage{subfigure}
\usepackage[utf8]{inputenc}
\newtheorem{obs}{Observation}
\newcommand{\Tt}[1]{\mathbf{#1}}

\begin{document}

\title{Wave propagation characteristics of Parareal
}


\author{Daniel Ruprecht 
}

\authorrunning{} 

\institute{Daniel Ruprecht \at
School of Mechanical Engineering, University of Leeds, LS2 9JT, UK
}

\date{Received: date / Accepted: date}

\maketitle

\begin{abstract}
The paper derives and analyses the (semi-)discrete dispersion relation of the Parareal parallel-in-time integration method.
It investigates Parareal's wave propagation characteristics with the aim to better understand what causes the well documented stability problems for hyperbolic equations.
The analysis shows that the instability is caused by convergence of the amplification factor to the exact value from above for medium to high wave numbers.
Phase errors in the coarse propagator are identified as the culprit, which suggests that specifically tailored coarse level methods could provide a remedy.
\keywords{Parareal \and convergence \and dispersion relation \and hyperbolic problem \and advection-dominated problem}
\end{abstract}



\section{Introduction}
Parallel computing has become ubiquitous in science and engineering, but requires suitable numerical algorithms to be efficient.
Parallel-in-time integration methods have been identified as a promising direction to increase the level of concurrency in computer simulations that involve the numerical solution of time dependent partial differential equations~\cite{DongarraEtAl2014}.
A variety of methods has been proposed~\cite{EmmettMinion2012,FalgoutEtAl2014_MGRIT,FarhatEtAl2003,Kiehl1994,LionsEtAl2001}, the earliest going back to 1964~\cite{Nievergelt1964}.
While even complex diffusive problems can be tackled successfully~\cite{FalgoutEtAl2016_sisc,FischerEtAl2005,KreienbuehlEtAl2015,Trindade2004} -- although parallel efficiencies remain low -- hyperbolic or advection-dominated problems have proved to be much har\-der to parallelise in time.
This currently prevents the use of pa\-rallel-in-time integration for most problems in computational fluid dynamics, even though many applications struggle with excessive solution times and could benefit great\-ly from new parallelisation strategies.

For the Parareal parallel-in-time algorithm there is some theory available illustrating its limitations in this respect.
Bal shows that Para\-real with a sufficiently damping coarse method is unconditionally stable for parabolic problems but not for hyperbolic equations~\cite{Bal2005}.
An early investigation of Parareal's stability properties showed instabilities for imaginary eigenvalues~\cite{StaffRonquist2005}.
Gander and Vandewalle~\cite{GanderVandewalle2007} give a detailed analysis of Parareal's convergence and show that even for the simple advection equation $u_t + u_x = 0$, Parareal is either unstable or inefficient.
Numerical experiments reveal that the instability em\-erges in the nonlinear case as a degradation of convergence with increasing Reynolds number~\cite{SteinerEtAl2015}.
Approaches exist to stabilise Parareal for hyperbolic equations~\cite{ChenEtAl2014,DaiEtAl2013,EghbalEtAl2016,FarhatCortial2006,GanderPetcu2008,KooijEtAl2017,RuprechtKrause2012}, but typically with significant overhead, leading to further degradation of efficiency, or limited applicability.

Since a key characteristic of hyperbolic problems is the existence of waves propagating with finite speeds, understanding Parareal's wave propagation characteristics is important to understand and, hopefully, resolve these problems.
However, no such analysis exists that gives insight into \emph{how} the instability emerges.
A better understanding of the instability could show the way to novel methods that allow the efficient and robust parallel-in-time solution of flows governed by advection.
Additionally, just like for ``classical'' time stepping methods, detailed knowledge of Parareal's theoretical properties for test cases will help understanding its performance for complex test problems where mathematical theory is not available.

To this end, the paper derives a discrete dispersion relation for Parareal to study how plane wave solutions $u(x,t) = \exp(-i \omega t) \exp(i \kappa x)$ are propagated in time.
It studies the discrete phase speed and amplification factor and how they depend on e.g. the number of processors, choice of coarse propagator and other parameters.
The analysis reveals that the source of the instability is a convergence from above in the amplification factor in higher wave number modes.
In diffusive problems, where high wave numbers are naturally strongly damped, this does not cause any problems, but in hyperbolic problems with little or no diffusion it causes the amplification factors to exceed a value of one and thus triggers instability.
Furthermore, the paper identifies phase errors in the coarse propagator as the source of these issues. 
This suggests that controlling coarse level phase errors could be key to devising efficient parallel-in-time methods for hyperbolic equations.

All results presented in this paper have been produced using \texttt{pyParareal}, a simple open-source Python code.
It is freely available~\cite{pyparareal} to maximise the usefulness of the here presented analysis, allowing other researchers to test different equations or to explore sets of parameters which are not analysed in this paper.


\section{Parareal for linear problems}
Parareal~\cite{LionsEtAl2001} is a parallel-in-time integration method for an initial value problem
\begin{equation}
	\label{eq:ivp}
	\dot{u}(t) = A u(t), \quad u(0) = u_0 \in \mathbb{C}^n, \quad 0 \leq t \leq T.
\end{equation}
For the sake of simplicity we consider only the linear case where the right hand side is given by a matrix $A \in \mathbb{C}^{n,n}$.
To parallelise integration of~\eqref{eq:ivp} in time, Parareal decomposes the time interval $[0,T]$ into $P$ time slices
\begin{equation}
	[0, T] = [0, T_1) \cup [T_1, T_2) \cup \ldots \cup [T_{P-1}, T ],
\end{equation}
with $P$ indicating the number of processors.
Denote as $\mathcal{F}_{\delta t}$ and $\mathcal{G}_{\Delta t}$ two ``classical'' time integration methods with time steps of length $\delta t$ and $\Delta t$ (e.g. Runge-Kutta methods).
For the sake of simplicity, assume that all slice $[T_{j-1}, T_{j})$ have the same length $\Delta T$ and that this length is an integer multiple of both $\delta t$ and $\Delta t$ so that $\Delta T = N_{\text{c}} \Delta t$ and $\Delta T = N_{\text{f}} \delta t$.
Below, $\delta t$ will always denote the time step size of the fine method and $\Delta t$ the time step size of the coarse method, so that we omit the indices and just write $\mathcal{G}$ and $\mathcal{F}$ to avoid clutter.
Standard serial time marching using the method denoted as $\mathcal{F}$ would correspond to evaluating
\begin{equation}
	\label{eq:time_stepping}
	u_{p} = \mathcal{F}(u_{p-1}), \quad p=1, \ldots, P,
\end{equation}
where $u_p \approx u(T_p)$.
Instead, after an initialisation procedure to provide values $u^0_p$ -- typically running the coarse method once -- Parareal computes the iteration
\begin{equation}
	\label{eq:parareal_iteration}
	u^{k}_{p} = \mathcal{G}(u^k_{p-1}) + \mathcal{F}(u^{k-1}_{p-1}) - \mathcal{G}(u^{k-1}_{p-1}), \quad p=1, \ldots, P
\end{equation}
for $k=1, \ldots, K$ where the computationally expensive evaluation of the fine method can be parallelised over $P$ processors.
If the number of iterations $K$ is small enough and the coarse method much cheaper than the fine, iteration~\eqref{eq:parareal_iteration} can run in less wall clock time than serially computing~\eqref{eq:time_stepping}.

\subsection{The Parareal iteration in matrix form}
As a first step toward deriving Parareal's dispersion relation we will need to derive its stability function which will require writing it in matrix form.
Consider now the case where both the coarse and the fine integrator are one-step methods with stability functions $R_f$ and $R_c$.
Then, $\mathcal{G}$ and $\mathcal{F}$ can be expressed as matrices 
\begin{equation}
	u_{p} = \mathcal{F}(u_{p-1}) = F u_{p-1}, \quad u_{p} = \mathcal{G}(u_{p-1}) = G u_{p-1}
\end{equation}
with $F := \left( R_f(A \delta t) \right)^{N_f}$ and $G := \left( R_c(A \Delta t) \right)^{N_c}$.
Denote as $\Tt{u}^k = ( u_0, \ldots, u_{P} ) \in \mathbb{R}^{(P+1)n}$ a vector that contains the approximate solutions at all time points $T_j$, $j=1, \ldots, P$ and the initial value $u_0$.
Simple algebra shows that one step of iteration~\eqref{eq:parareal_iteration} is equivalent to the block matrix formulation
\begin{equation}
	\label{eq:parareal_iteration_matrix}
	\Tt{M}_{g} \Tt{u}^{k} = \left( \Tt{M}_{g} - \Tt{M}_{f} \right) \Tt{u}^{k-1} + \Tt{b}	
\end{equation}
with matrices
\begin{equation}
	\Tt{M}_{f} := \begin{bmatrix} I & & & \\ -F & I & \\  & \ddots & \ddots \\ & & -F & I \end{bmatrix} \in \mathbb{R}^{(P+1)n, (P+1)n}
\end{equation}
and
\begin{equation}
	\Tt{M}_{g} := \begin{bmatrix} I & & & \\ -G & I & \\  & \ddots & \ddots \\ & & -G & I \end{bmatrix} \in \mathbb{R}^{(P+1)n, (P+1)n}
\end{equation}
and a vector $\Tt{b} = ( u_0, 0, \ldots, 0 ) \in \mathbb{R}^{(P+1)n}$.
Formulation~\eqref{eq:parareal_iteration_matrix} interpretes Parareal as a preconditioned linear iteration~\cite{AmodioBrugnano2009}.

\subsection{Stability function of Parareal}
From the matrix formulation of a single Parareal iteration~\eqref{eq:parareal_iteration_matrix}, we can now derive its stability function, that is we can express the update from the initial value $u_0$ to an approximation $u_P$ at time $T=T_P$ using Parareal with $K$ iterations as multiplication by a single matrix.
The fine propagator solution satisfies
\begin{equation}
	\label{eq:fine_solution_equation}
	\Tt{M}_f \Tt{u}_f = \Tt{b}
\end{equation}
and is a fixed point of iteration~\eqref{eq:parareal_iteration_matrix}.
Therefore, propagation of the error
\begin{equation}
	\Tt{e}^k := \Tt{u}^k - \Tt{u}_f
\end{equation}
is governed by the matrix
\begin{equation}
	\label{eq:error_matrix}
	\Tt{E} := \left( \Tt{I} - \Tt{M}_g^{-1} \Tt{M}_{f} \right)
\end{equation}
in the sense that
\begin{equation}
	\label{eq:error_propagation}
	\Tt{e}^{k} = \Tt{E} \Tt{e}^{k-1}.
\end{equation}
Using this notation and applying~\eqref{eq:parareal_iteration_matrix} recursively, it is now easy to show that
\begin{equation}
	\label{eq:closed_form_1}
	\Tt{u}^k  =  \Tt{E} \Tt{u}^0 + \sum_{j=0}^{k-1} \Tt{E}^j \Tt{M}_{g}^{-1} \Tt{b}.
\end{equation}
If, as is typically done, the start value $\Tt{u}^0$ for the iteration is produced by a single run of the coarse method, that is if
\begin{equation}
	\label{eq:coarse_solution_equation}
	\Tt{M}_g \Tt{u}^0 = \Tt{b},
\end{equation}
Equation~\eqref{eq:closed_form_1} further simplifies to
\begin{equation}
	\label{eq:closed_form_2}
	\Tt{u}^k = \left( \sum_{j=0}^{k} \Tt{E}^j \right) \Tt{M}_g^{-1} \Tt{b}.
\end{equation}
The right hand side vector can be generated from the initial value $u_0$ via
\begin{equation}
	\Tt{b} = \Tt{C}_1 u_0
\end{equation}
by defining
\begin{equation}
	\Tt{C}_1 = \left[ \Tt{I} ; \Tt{0} \right] \in \mathbb{R}^{(P+1)n,n}, \quad \Tt{I} \in \mathbb{R}^{n,n}, \Tt{0} \in \mathbb{R}^{Pn, n}.
\end{equation}
Finally, denote as 
\begin{equation}
	\Tt{C}_2 = \left[ \Tt{0}, \Tt{I} \right], \quad \Tt{0} \in \mathbb{R}^{n, Pn}, \ \Tt{I} \in \mathbb{R}^{n,n},
\end{equation}
the matrix that selects the last $n$ entries out of $\Tt{u}^k$.
Now, a full Parareal update from some initial value $u_0$ to an approximation $u_{P}$ using $K$ iterations can be written compactly as
\begin{equation}
	\label{eq:parareal_update}
	u_{P} = \Tt{C}_2 \left( \sum_{j=0}^{K} \Tt{E}^j \right) \Tt{M}_g^{-1} \Tt{C}_1 u_0 =: M_{\text{parareal}} u_0.
\end{equation}
The stability matrix $M_{\text{parareal}} \in \mathbb{R}^{n,n}$ depends on $K$, $T$, $\Delta T$, $P$, $\Delta t$, $\delta t$, $\mathcal{F}$,  $\mathcal{G}$ and $A$.
Note that Staff and R{\o}nquist derived the stability function for the scalar case using Pascal's tree~\cite{StaffRonquist2005}.

\subsection{Weak scaling vs. longer simulation times}
There are two different application scenarios for Parareal that we can study when increasing the number of processors $P$.
If we fix the final time $T$, increasing $P$ will lead to better resolution since the coarse time step $\Delta t$ cannot be larger than the length of a time slice $\Delta T$ -- the coarse method has to perform at least one step per time slice.
In this scenario, more processors are used to absorb the cost of higher temporal resolution (``weak scaling'').

Alternatively, we can use additional processors to compute until a later final time $T$ and this is the scenario investigated here.
Consequently, we study here the case where $T$ and $P$ increase together and always assume that $T=P$, that is each time slice has length one and increasing $P$ means parallelising over more time slices covering a longer time interval.
Since dispersion properties of numerical methods are typically analysed for a unit interval, this causes some issues that we resolve by ``normalising'' the Parareal stability function, see \S\ref{subsec:normalisation}.

\subsection{Maximum singular value}
\begin{figure}[t!]
	\centering
	\includegraphics[scale=1]{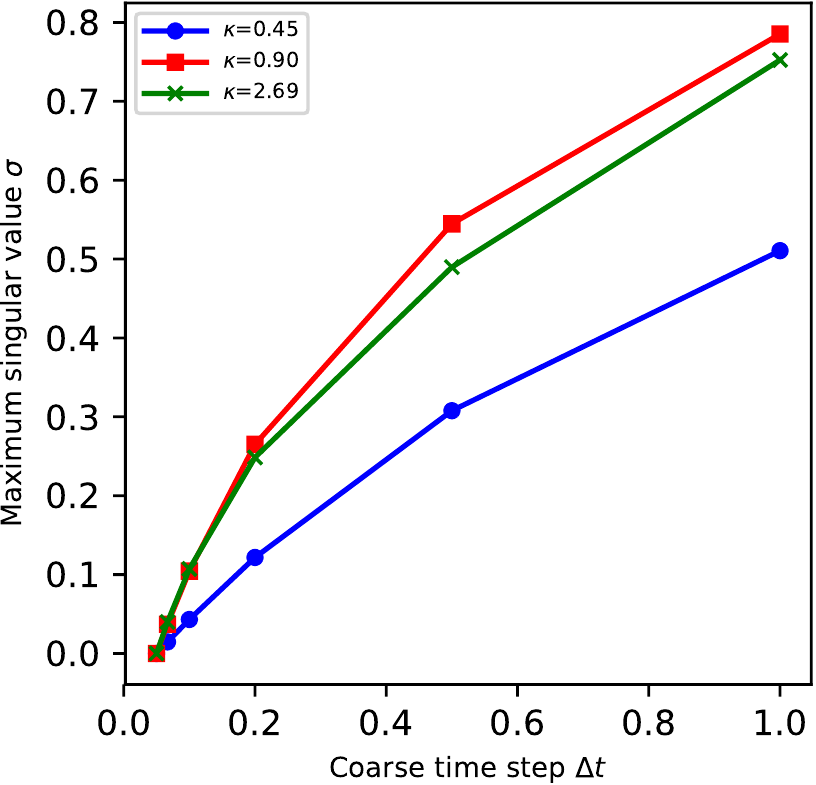}
	\caption{Maximum singular value $\sigma$ for decreasing coarse time step $\Delta t$ for $\nu = 0$. Backward Euler is used for both $\mathcal{F}$ and $\mathcal{G}$ and for $\Delta t = \delta t = 0.05$, both methods coincide so that $\sigma = 0$.}
	\label{fig:svd-vs-dt}
\end{figure}

The matrix $\Tt{E}$ defined in~\eqref{eq:error_matrix} determines how quickly Parareal converges.
Note that $\Tt{E}$ is nil-potent with $\Tt{E}^{P} = 0$, owing to the fact that after $P$ iterations Parareal will always reproduce the fine solution exactly.
Therefore, all eigenvalues of $\Tt{E}$ are zero and the spectral radius is not useful to analyse convergence.
Below, to investigate convergence, we will therefore compute the maximum singular value $\sigma$ of $\Tt{E}$ instead.
Since
\begin{equation*}
	\sigma = \left\| \Tt{E} \right\|_2,
\end{equation*}
it follows from~\eqref{eq:error_propagation} that
\begin{equation}
	\label{eq:svd_bound}
	\left\| \Tt{e}^{k} \right\|_2 \leq \left\| \Tt{E} \right\|_2 \left\| \Tt{e}^{k-1} \right\|_2 = \sigma \left\| \Tt{e}^{k-1} \right\|_2 \leq \sigma^{k} \left\| \Tt{e}^{0} \right\|_2
\end{equation}
so that if $\sigma < 1$ Parareal will converge monotonically without stal\-ling.
In particular, this rules out behaviour as found by Gander and Vandewalle for hyperbolic problems, where the error would first increase substantially over the first $P/2$ iterations before beginning to decrease~\cite{GanderVandewalle2007_SISC}.
However, achieving fast convergence and good efficiency will typically require $\sigma \ll 1$.
Note that if the coarse method is used to generate $\Tt{u}_0$, it follows from~\eqref{eq:fine_solution_equation} and~\eqref{eq:coarse_solution_equation} that the initial error is
\begin{equation}
	\Tt{e}^0 = \Tt{u}^{0} - \Tt{u}_f = \left( \Tt{M}_g^{-1} - \Tt{M}_f^{-1} \right) \Tt{b}.
\end{equation}
The size of $\sigma$ depends on the accuracy ``gap'' between the coarse and fine integrator and the wave number.
Figure~\ref{fig:svd-vs-dt} shows $\sigma$ for varying values of $\Delta t$ when backward Euler is used for both coarse and fine method.
Clearly, as the coarse time step approaches to fine time step of $\delta t = 0.05$, the maximum singular value approaches zero.
However, in this limit the coarse and fine propagator are identical and no speedup is possible.
The larger $\Delta t$ compared to $\delta t$, the cheaper the coarse method becomes but since $\sigma$ also grows, more iterations are likely required.
Note that higher wave numbers lead to higher values of $\sigma$ while lower wave numbers tend to have values of $\sigma \ll 1$ even for large coarse-to-fine time step ratios.

\begin{figure}[t!]
	\centering
	\includegraphics[scale=1]{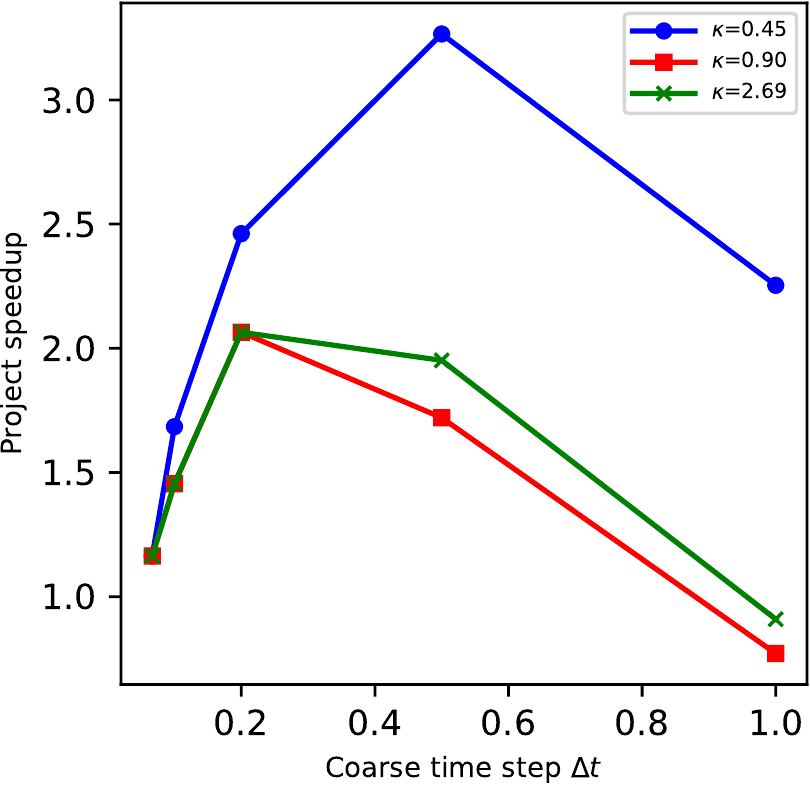}
	\caption{Projected speedup for pipelined Parareal~\cite{Minion2010} with $P=16$ processors for the same parameters as in Figure~\ref{fig:svd-vs-dt} and the number of iterations $k$ fixed such that $\sigma^k \leq \text{tol}$ with $\text{tol} = 0.01$.}	
	\label{fig:speedup-vs-dt}
\end{figure}

Looking at $\sigma$ also provides a way to refine performance models for Parareal.
Typically, in models projecting speedup, the number of iterations has to be fixed in addition to $\Delta t$, $\delta t$ and $P$.
Instead, at least for linear problems, we can fix $k$ such that
\begin{equation}
	\sigma^{k} \leq \text{tol}
\end{equation}
for some fixed tolerance $\text{tol}$.
The resulting projected speedup for $P=16$ processors and a tolerance of $\text{tol} = 1e-2$ is shown in Figure~\ref{fig:speedup-vs-dt}.
First, as the coarse time step increases, the reduced cost of the coarse propagator improves achievable speedup.
Simultaneously, the decreasing accuracy of $\mathcal{G}$ manifests itself in an increasing number of iterations required to match the selected tolerance.
These two counteracting effects create a ``sweet spot`` where $\mathcal{G}$ is accurate enough to still enable relatively fast convergence but cheap enough to allow for speedup.
It is noteworthy, however, that this sweet spot is different for lower and higher wave numbers.
Therefore, the potential for speedup from Parareal does not solely depend on the solved equations and discretization parameters but also on the solution - the more prominent high wave number modes are, the more restricted achievable speed\-up.

\subsection{Convergence and (in)stability of Parareal}
Two different but connected issues with Parareal are discussed throughout the paper, convergence and (in)stability.
Here, convergence refers to how fast Parareal approaches the fine solution within a single instance of Parareal, that is
\begin{equation}
	M_{\text{parareal}} \to F \quad \text{as} \quad k \to P.
\end{equation}
As discussed above, after $k=P$ iterations we always have $M_{\text{parareal}} = F$, but particularly for hyperbolic problems this can happen only at the final iterations $k=P$ while at $k=P-1$ there is still a substantial difference~\cite{GanderVandewalle2007_SISC}.
Clearly, speedup is not obtainable in such a situation.
The maximum singular value $\sigma$ of $\Tt{E}$ gives an upper bound or worst-case scenario of how fast Parareal converges to the fine solution, cf. Equation~\eqref{eq:svd_bound}.
While $\sigma < 1$ does not necessarily guarantee converge quick enough to generate meaningful speedup, it guarantees monotonous convergence and rules out an error that increases first before decreasing only in later iterations.

The other issue investigated in the paper is that of stability of repeated application of Parareal (``restarting'').
Below, stability is normally assessed for Parareal with a fixed number of iterations $k$.
A configuration of Parareal is referred to as \emph{unstable} if it leads to an amplification factor of more than unity.
This corresponds to an artificial increase in wave amplitudes and, just as for classical time stepping methods, would result in a diverging numerical approximation if the method is used recursively
\begin{equation}
	\left( M_{\text{parareal}} \right)^n \to \infty \quad \text{as} \quad n \to \infty.
\end{equation}
While for classical methods this recursive application simply means stepping through time steps, for Parareal with $P$ processors it would mean computing one window $[0, T_{P}]$, then restarting it with the final approximation as initial value for the next window $[T_P, 2 T_{P}]$ and so on.
\section{Discrete dispersion relation for the advection-diffusion equation}
Starting from~\eqref{eq:parareal_update} we now derive the (semi-)discrete dispersion relation of Parareal for the one dimensional linear advection diffusion problem
\begin{equation}
	\label{eq:advection_diffusion}
	u_t + U u_x = \nu u_{xx}.
\end{equation}
First, assume a plane wave solution in space
\begin{equation}
	\label{eq:cont_plane_wave}
	u(x,t) = \hat{u}(t) e^{i \kappa x}
\end{equation}
with wave number $\kappa$ so that~\eqref{eq:advection_diffusion} reduces to the initial value problem
\begin{equation}
	\label{eq:u_hat}
	\hat{u}_t(t) = -\left( U i \kappa + \nu \kappa^2 \right) \hat{u}(t) =: \delta(\kappa, U, \nu) \hat{u}(t) 
\end{equation}
with initial value $\hat{u}(0) = 1$.
Integrating~\eqref{eq:u_hat} from $t=0$ to $t=T$ in one step gives
\begin{equation}
	\label{eq:u_hat_stabfunc}
	\hat{u}_{T} = R(\delta, T) \hat{u}_0
\end{equation}
where $R$ is the stability function of the employed method.
Now assume that the approximation of $\hat{u}$ is a discrete plane wave so that the solution at the end of the $n$\textsuperscript{th} time slice is given by
\begin{equation}
	\hat{u}_n = e^{-i \omega n \Delta T}.
\end{equation}
Inserting this in~\eqref{eq:u_hat_stabfunc} gives
\begin{equation}
	\label{eq:omega_log}
	e^{-i \omega T} \hat{u}_0 = R \hat{u}_0 \Rightarrow \omega = i \frac{\log(R)}{T}.
\end{equation}
For $R$ in polar coordinates, that is $R = \left| R \right| \exp(i \theta)$ with $\theta = \texttt{angle}(R)$, we get
\begin{equation}
	\omega = i \left( \log(\left| R \right| ) + i \theta \right) T^{-1}.	
\end{equation}
The exact integrator, for example, would read 
\begin{equation}
	\label{eq:exact_integrator}
	R_{\text{exact}} = e^{\delta(\kappa,U,\nu) T}.
\end{equation}
Using~\eqref{eq:omega_log} to compute the resulting frequency yields $\omega =  i \delta(k, U, \nu)$ and retrieves the continuous plane wave solution
\begin{equation}
	u(x,T) = \hat{u}_T e^{i \kappa x} =  e^{-\nu \kappa^2 T} e^{i \kappa \left( x - U T \right)}
\end{equation}
 of~\eqref{eq:advection_diffusion}.
It also reproduces the dispersion relation of the continuous system
\begin{equation}
	\label{eq:exact_dispersion_relation}
	\omega = i \frac{\log(R)}{T} = i  \delta(\kappa, U, \nu) =  U \kappa - i \nu \kappa^2.
\end{equation}
However, if we use an approximate stability function $R$ instead, we get some approximate $\omega = \omega_{r} + i \omega_{i}$ with $\omega_{r}, \omega_{i} \in \mathbb{R}$. 
The resulting semi-discrete solution becomes
\begin{equation}
	u_j^n = e^{-i \omega t_n} e^{i \kappa x} = e^{\omega_i t_n} e^{i \kappa \left( x - \frac{\omega_{r}}{\kappa} t_n \right)}.
\end{equation}
Therefore, $\omega_{r}/\kappa$ governs the propagation speed of the solution while $\omega_{i}$ governs the growth or decay in amplitude.
Consequently, $\omega_{r}/\kappa$ is referred to as \emph{phase velocity} while $\exp(\omega_i)$ is called \emph{amplification factor}.
In the continuous case, the phase speed is equal to $U$ and the amplification factor equal to $e^{-\nu \kappa^2}$.
Note that for~\eqref{eq:advection_diffusion} the exact phase speed should be independent of the wave number $\kappa$.
However, the discrete phase speed of a numerical method often will change with $\kappa$, thus introducing numerical dispersion.
Also note that for $\nu > 0$ higher wave numbers decay faster because the amplification factor decreases rapidly as $\kappa$ increases.

\subsection{Normalisation}\label{subsec:normalisation}
The update function R for Parareal in Equation~\eqref{eq:parareal_update} denotes not an update over $[0,1]$ but over $[0,T_P]$ where $T_P = P$ is the number of processors.
A phase speed of $\omega_{r}/\kappa = 1.0$, for example, indicates a wave that travels across a whole interval $[0,1]$ during the step.
If scaled up to an interval $[0,P]$, the corresponding phase speed would become $\omega_{r}/\kappa = P$ instead.

This enlarged range of values causes problems with the complex logarithm in~\eqref{eq:omega_log}.
As an example, take the stability function of the exact propagator~\eqref{eq:exact_integrator}.
Analytically, the identity
\begin{equation}
	\label{eq:dispersion_identity}
	\omega = i \frac{R}{T} = i \delta(\kappa, U, \nu) \frac{T}{T} = i \delta(\kappa, U, \nu)
\end{equation}
holds, resulting in the correct dispersion relation~\eqref{eq:exact_dispersion_relation} of the continuous system.
However, depending on the values of $\kappa$, $U$, $T$ and $\nu$, this identity is not necessarily satisfied when computing the complex logarithm using \texttt{np.log}.
For example, for $\kappa = 2$, $U = 1$,  $\nu = 0$ and $T>0$, the exact stability function is $R = e^{-2 i T}$.
In Python, we obtain for $T=1$
\begin{equation}
	\texttt{1j} * \texttt{np.log}(e^{-2 i})/1 = 2 = \kappa
\end{equation}
but for $T=2$ we obtain
\begin{equation}
	\texttt{1j} * \texttt{np.log}(e^{-4 i})/2 \approx -1.1416 \neq \kappa
\end{equation}
and so identity~\eqref{eq:dispersion_identity} is not fulfilled.
The reason is that the logarithm of a complex number is not unique and so \texttt{np.log} returns only \emph{a} complex logarithm but not necessarily the right one in terms of the dispersion relation.

To circumvent this problem, we ``normalise'' the update $R$ for Parareal to $[0,1]$.
To this end, decompose 
\begin{equation}
	R = \sqrt[P]{R} \cdot \ldots \cdot \sqrt[P]{R}
\end{equation}
where $\sqrt[P]{R}$ corresponds to the propagation over $[0,1]$ instead of $[0,P]$.
Since there are $P$ many roots $\sqrt[P]{R}$, we have to select the right one.
First, we use the \texttt{zeros} function of \texttt{numPy} to find all $P$ complex roots $z_i$ of
\begin{equation}
	z^n - R = 0.
\end{equation}
Then, we select as root $\sqrt[P]{R}$ the value $z_i$ that satisfies
\begin{equation}
	 \left| \theta(\sqrt[P]{R}) - \theta_{\rm targ} \right| = \min_{p=1,\ldots,P} \left| \theta(z_p) - \theta_{\rm targ} \right|
\end{equation} 
where $\theta$ is the \texttt{angle} function and $\theta_{\rm targ}$ some target angle, which we still need to define.

We compute $\omega$ and the resulting phase speed and amplification factor for a range of wave numbers $0 \leq \kappa_1 \leq \kappa_2 \leq \ldots \leq \kappa_{N} \leq \pi$.
For $\kappa_1$, $\theta_{\rm targ}$ is set to the angle of the frequency $\omega$ computed from the analytic dispersion relation.
After that, $\theta_{\rm targ}$ is set to the angle of the root selected for the previous value of $\kappa$.
The rationale is that small changes in $\kappa$ should only result in small changes of frequency and phase so that $\theta(\omega_{i-1}) \approx \theta(\omega_i)$ if the increment between wave numbers is small enough.
From the selected root $\sqrt[P]{R}$ we then compute $\omega$ using~\eqref{eq:omega_log}, the resulting discrete phase speed and amplification factor and the target angle $\theta_{\rm targ}$ for the next wave number.

\section{Analysis of the dispersion relation}
After showing how to derive Para\-real's dispersion relation and normalising it to the unit interval, this section now provides a detailed analysis of different factors influencing its discrete phase speed and amplification factor.

\subsection{Influence of diffusion}
\begin{figure*}[ht!]
	\centering
	\includegraphics[scale=1]{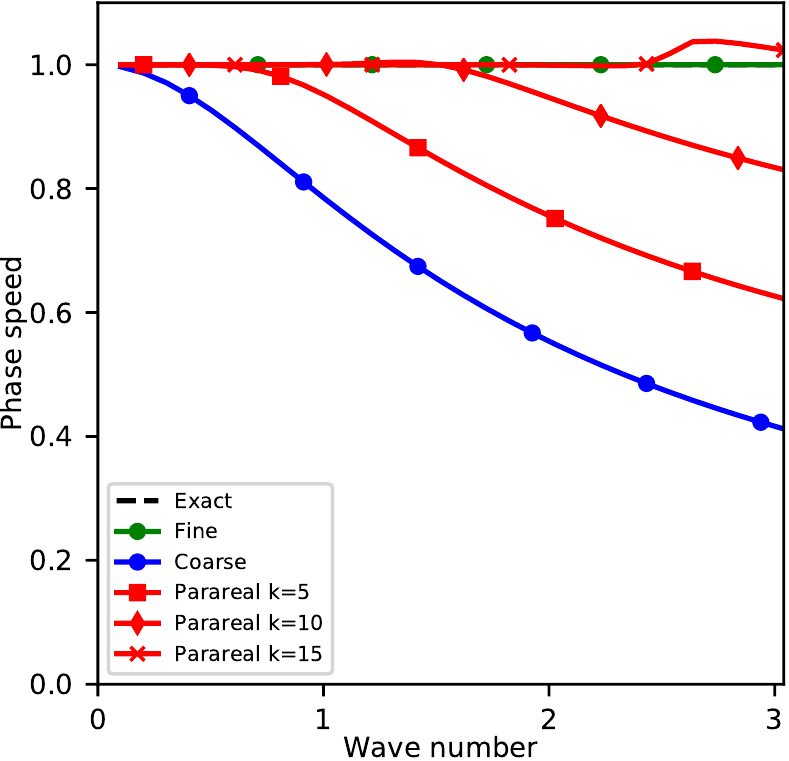}\hfill
	\includegraphics[scale=1]{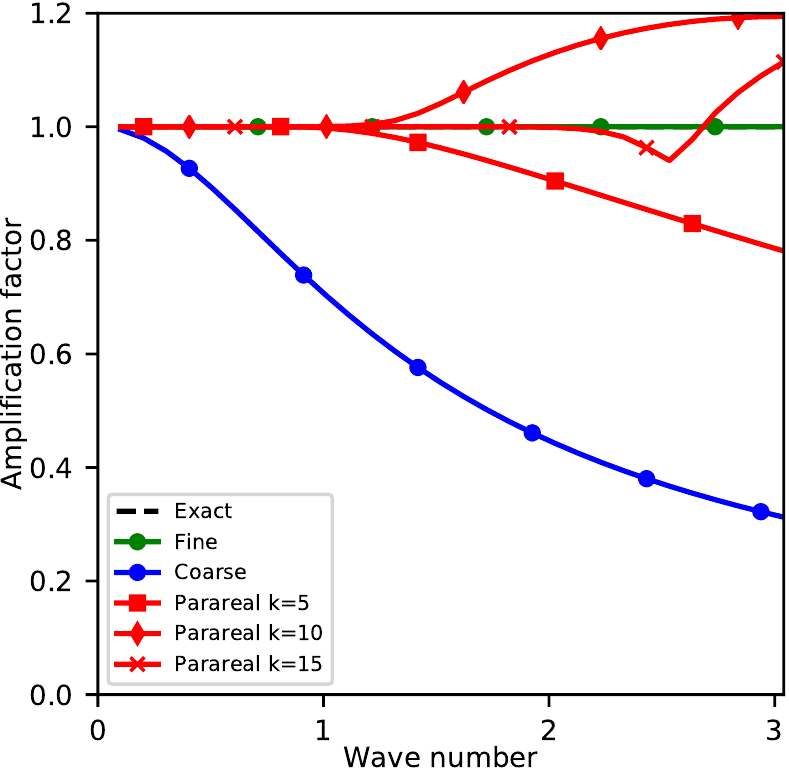}
	\includegraphics[scale=1]{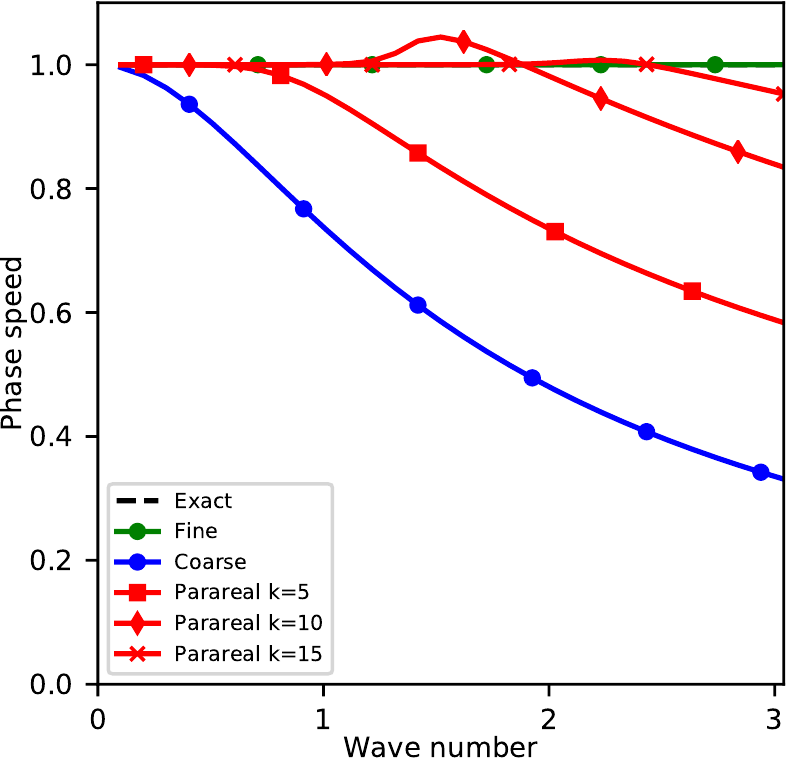}\hfill
	\includegraphics[scale=1]{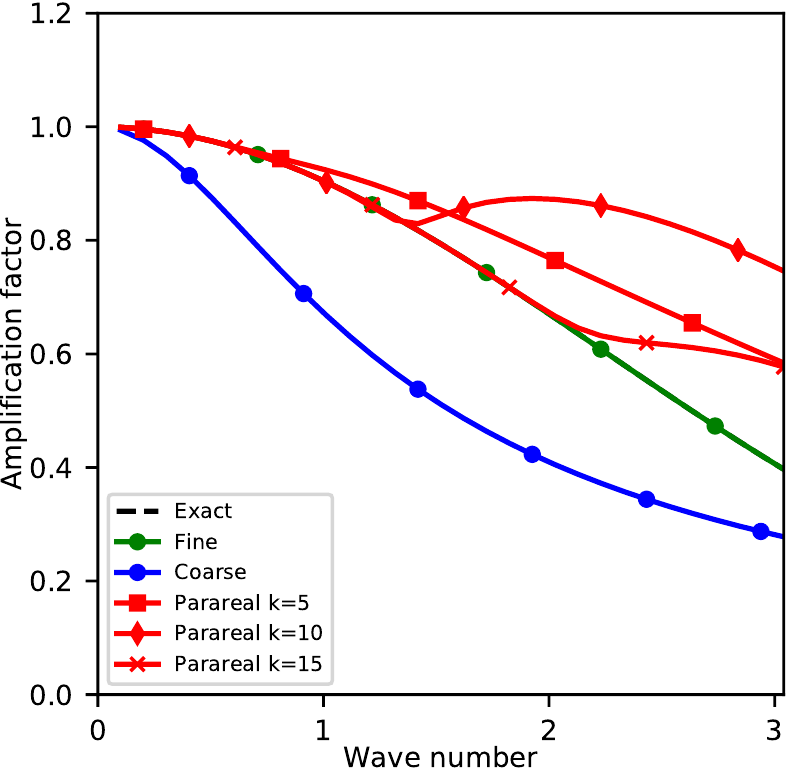}
	\caption{Discrete phase speed and amplification factor for Parareal with backward Euler as $\mathcal{G}$ and the exact integrator for $\mathcal{F}$. The symbol for the spatial discretisation is $\delta = - \left( i \kappa + \nu \kappa^2\right)$. The diffusion coefficient is $\nu=0.0$ (upper) and $\nu=0.1$ (lower).}
	\label{fig:dispersion}
\end{figure*}
Figure~\ref{fig:dispersion} shows the discrete phase speed and amplification factor of Parareal for $P=16$, backward Euler with $\Delta t = 1.0$ as coarse and the exact integrator as fine propagator.
Both levels use $\delta = -\left( U i \kappa + \nu \kappa^2 \right)$, that is the symbol of the continuous spatial operator.
The two upper figures are for $U=1.0$ and $\nu=0.0$ (no diffusion) while the two lower figures are for $U=1.0$ and $\nu=0.1$ (diffusive).

In both cases, the discrete phase speed of Parareal converges almost monotonically toward the continuous phase speed.
Even for ten iterations, Parareal still causes significant slowing of medium to large wave number modes.
Parareal requires almost the full number of iterations, $k=15$, before it faithfully reproduces the correct phase speed across most of the spectrum.
However, for any number of iterations where speed up might still be possible, Parareal will introduce significant numerical dispersion.
Slight artificial acceleration is also observed for high wave number modes for $k=15$ in the non-diffusive and $k=10$ in the diffusive case, but generally phase speeds are quite similar in the diffusive and non-diffusive case.

The amplification factor in the non-diffusive case (upper right figure) illustrates Parareal's instability for hyperbolic equations: for $k=10$ and $k=15$ it is larger than one for a significant part of the spectrum, indicating that these modes are unstable and will be artificially amplified.
For $k=5$, the iteration has not yet corrected for the strong diffusivity of the coarse propagator and remains stable for all modes but with significant numerical damping of medium to high wave numbers.
The reason for the stability problems is discernible from the amplification factor for the diffusive case (rightmost figure): from $k=0$ (blue circles) to $k=5$, Parareal reproduces the correct amplification factor for small wave number modes but significantly overestimates the amplitude of medium to large wave numbers.
It then continues to converge to the correct value \emph{from above}.
For the diffusive case where the exact values are smaller than one this does not cause instabilities. 
In the non-diffusive case, however, any overestimation of the analytical amplification factor immediately causes instability.
There is, in a sense, ``no room'' for the amplification factor to converge to the correct value from above.
This also means that using a non-diffusive method as coarse propagator, for example trapezoidal rule, leads to disastrous consequences (not shown) where most parts of the spectrum are unstable for almost any value of $k$.

\begin{figure*}[ht!]
	\centering
	\includegraphics[scale=1]{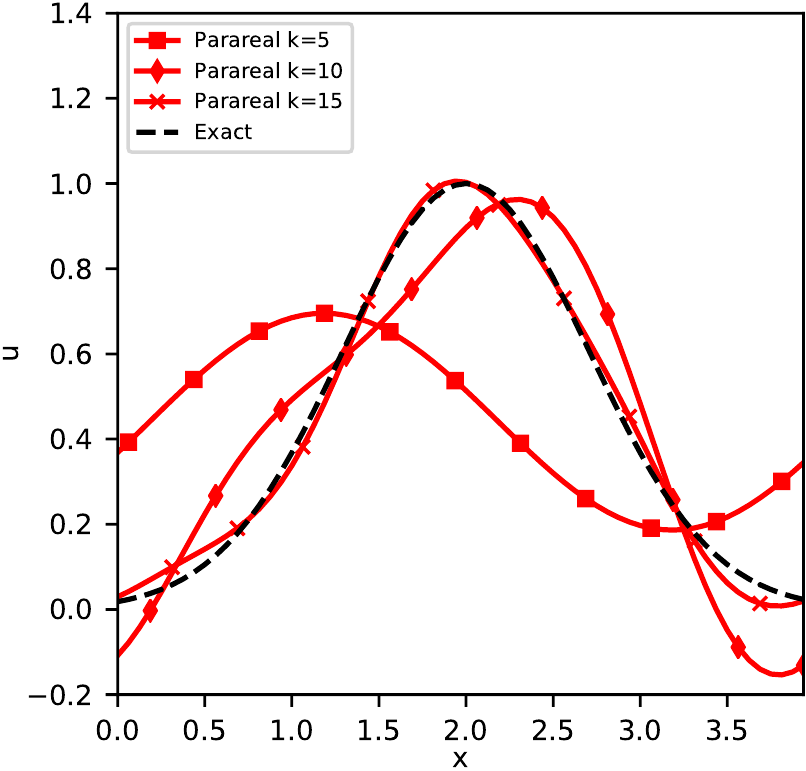}\hfill
	\includegraphics[scale=1]{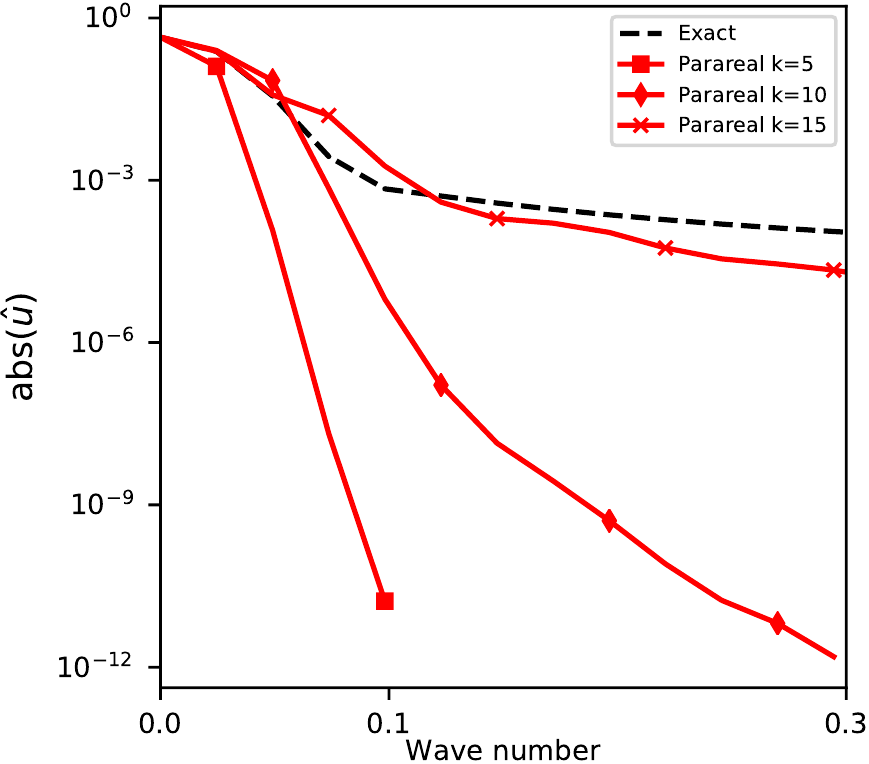}
	\caption{Gauss peak in physical space (left) and corresponding spectrum (right) for $U=1.0$ and $\nu=0.0$ integrated using a pseudo-spectral method with $64$ modes in space and Parareal with $P=16$ processors in time.}
	\label{fig:gauss_peak}
\end{figure*}

Figure~\ref{fig:gauss_peak} illustrate how the phase speed and amplitude errors discussed above manifest themselves.
It shows a single Gauss peak advected with a velocity of $U=1.0$ with $\nu = 0.0$ on a spatial domain $[0,4]$ over a time interval $[0,16]$ distributed over $P=16$ processors and $N_c = 2$ coarse time steps per slice.
A spectral discretisation is used in space, allowing to represent the derivative exactly.
For $k=5$ iterations, most higher wave numbers are damped out and the result looks essentially like a low wave number sine function.
The artificially amplified medium to high wave number modes create a ``bulge'' for $k=10$ while dispersion leads to a significant trough at the sides of the domain.
After fifteen iterations, the solution approximates the main part of the Gauss peak reasonably well, but dispersion still leads to visible errors along the flanks.
The right figure shows a part of the resulting spectrum.
For $k=5$, only the lowest wave number modes are present, leading to the sine shaped solution.
After $k=10$ iterations, most of the spectrum is still being truncated but a small range of wave numbers around $\kappa = 0.05$ is being artificially amplified which creates the ``bulge'' seen in the left figure.
Finally, for $k=15$ iterations, Parareal starts to correctly capture the spectrum but the still significant overestimation of  low wave number amplitudes and underestimation of higher modes causes visible errors.

\begin{obs}
The amplification factor in Parareal for high\-er wave numbers converges ``from above''. In diffusive problems these wave numbers are damped, so the exact amplification factor is significantly smaller than one, leaving room for Parareal to overestimate it without crossing the threshold to instability. For non-diffusive problems where the exact amplification factor is one, every overestimation causes the mode to become unstable.
\end{obs}

\subsection{Low order finite difference in coarse method}
In a realistic scenario, some approximation of the spatial derivatives would have to be used instead of the exact symbol $\delta$.
For simple finite differences, we can study the effect this has on the dispersion relation.
Consider the first order upwind finite difference
\begin{equation*}
	u_x(x_j) \approx \frac{u_j - u_{j-1}}{\Delta x}
\end{equation*}	
as approximation for $u_x$ in~\eqref{eq:advection_diffusion}.
Assuming a discrete plane wave 
\begin{equation*}
	u_j = \hat{u}(t) e^{i \kappa j \Delta x}
\end{equation*}
on a uniform spatial mesh $x_j = j \Delta x$ instead of the continuous plane wave~\eqref{eq:cont_plane_wave}, this leads to
\begin{equation*}
	u_x(x_j) \approx \frac{e^{i \kappa x_j} - e^{i \kappa (x_j - \Delta x)}}{\Delta x} = e^{i \kappa x_j} \frac{1 - e^{- i \kappa \Delta x}}{\Delta x}.
\end{equation*}
For $\nu=0$ this results in the initial value problem
\begin{equation*}
	\hat{u}_t(t) = - U \frac{1 - e^{- i \kappa \Delta x}}{\Delta x} \hat{u}(t) =: \tilde{\delta}(k, U, \delta x) \hat{u}(t)
\end{equation*}
with initial value $\hat{u}(0) = 1$ and a discrete symbol $\tilde{\delta}$ instead of $\delta$ as in~\eqref{eq:u_hat}.
Note that
\begin{equation*}
	\lim_{\Delta x \to 0} \frac{1 - e^{- i \kappa \Delta x}}{\Delta x} = i \kappa
\end{equation*}
so that $\tilde{\delta} \to \delta$ as $\Delta x \to 0$.
Durran gives details for different stencils~\cite{Durran2010}.

The dispersion properties of the implicit Euler method together with the first order upwind finite difference are qualitatively similar to the ones for implicit Euler with the exact symbol (not shown).\footnote{The script \texttt{plot\_ieuler\_dispersion.py} supplied together with the Python code can be used to visualize the dispersion properties of the coarse propagator alone.}
Using the upwind finite difference instead of the exact symbol gives qualitatively similar wave propagation characteristics for the coarse propagator.
Numerical slowdown increases (up to the point where modes at the higher end of the spectrum almost do not propagate at all) and numerical diffusion becomes somewhat stronger.
As a result, Parareal's dispersion properties (also not shown) are also relatively similar, except for too small phase speeds even for $k=15$.

\begin{figure*}[ht!]
	\centering
	\includegraphics[scale=1]{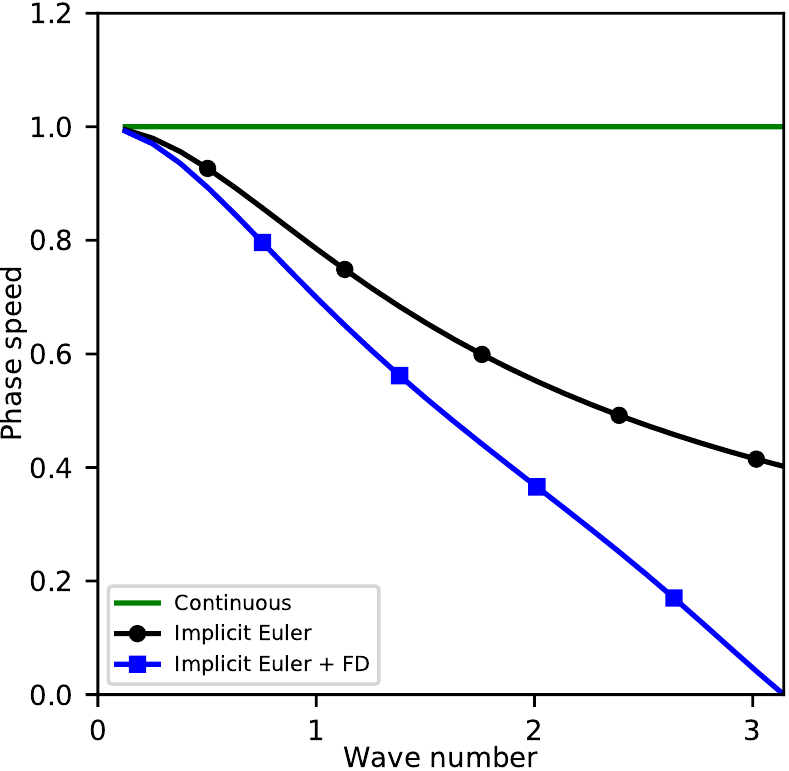}\hfill
	\includegraphics[scale=1]{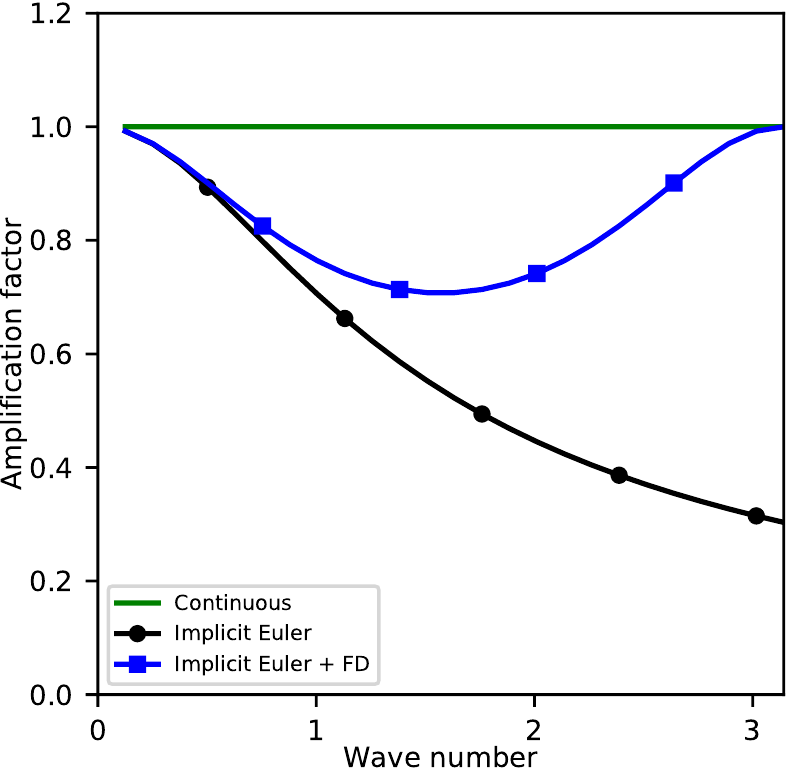}
	\caption{Phase speed (upper) and amplification factor (lower) of the implicit Euler method using the exact symbol $\delta$ (black circles) or the approximate symbol $\tilde{\delta}$ of the second order centred finite difference (blue squares).}
	\label{fig:dispersion-implicit-euler}
\end{figure*}

However, if we use the centred finite difference 
\begin{equation*}
	u_x(x_j) \approx \frac{u_{j+1} - u_{j-1}}{2 \Delta x}
\end{equation*}
instead, this leads to an approximate symbol
\begin{equation}
	\tilde{\delta} = -U \frac{e^{i \kappa \Delta x} - e^{-i\kappa \Delta x}}{2 \Delta x} = -U i \frac{\sin(\kappa \Delta x)}{\Delta x}.
\end{equation}
In this case, it turns out that the dispersion properties of implicit Euler with $\delta$ and $\tilde{\delta}$ are quite different.
Figure~\ref{fig:dispersion-implicit-euler} shows the discrete phase speed (left) and amplification factor (right) for both configurations.
For the phase speed, both version agree qualitatively, even though using $\tilde{\delta}$ leads to noticeable stronger slow down, particularly of higher wave numbers.
For the amplification factor, however, there is a significant difference between the semi-discrete and fully discrete method.
While the former damps high wave numbers strongly, the combination of implicit Euler and centred finite differences strongly damps medium wave numbers while damping of high wave numbers is weak.

\begin{figure*}[ht!]
	\centering
	\includegraphics[scale=1]{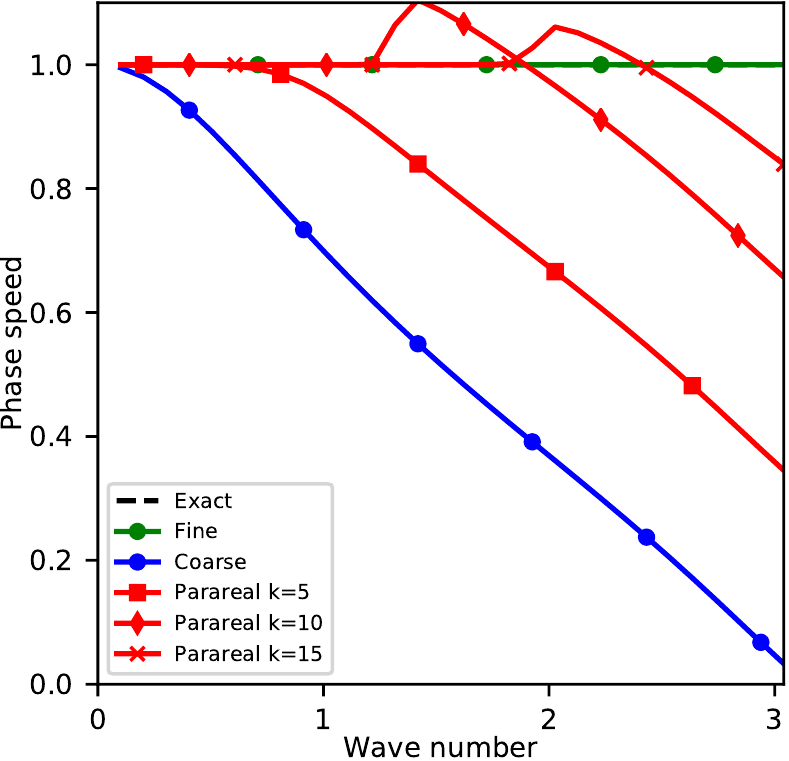}\hfill
	\includegraphics[scale=1]{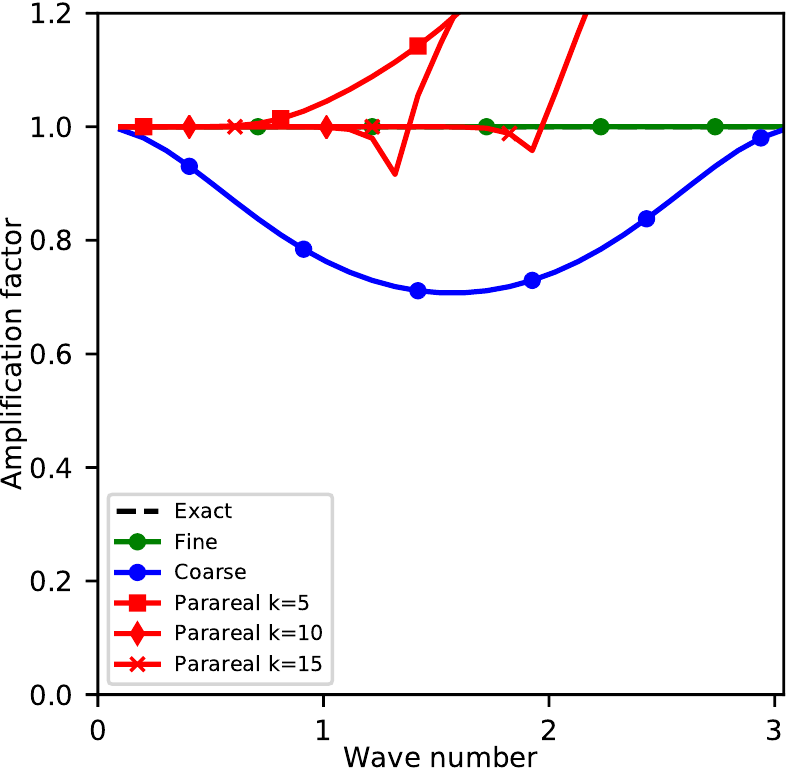}
	\caption{Phase speed (upper) and amplification factor (lower) for same configuration as in Figure~\ref{fig:dispersion} but with a second order centred finite difference in the coarse propagator instead of the exact symbol.}
	\label{fig:dispersion-fd-g}
\end{figure*}

In Parareal, this causes a situation similar to what happens when using the trapezoidal rule as coarse propagator, albeit less drastic.
Figure~\ref{fig:dispersion-fd-g} shows again the phase speed (left) and amplification factor (right) for the same configuration as before but implicit Euler with centred finite difference for $\mathcal{G}$.
The failure of the coarse method to remove high wave number modes again leads to an earlier triggering of the instability.
Whereas for Parareal using $\delta$ on the coarse level the iteration $k=5$ was will stable (see Figure~\ref{fig:dispersion}), it is now unstable.
For iterations $k=10$ and $k=15$, large parts of the spectrum remain unstable.
Also, the stronger numerical slow down of the coarse method makes it harder for Parareal to correct for phase speed errors.
Where before Parareal with $k=15$ iteration captured the exact phase speed reasonably well, in Figure~\ref{fig:dispersion-fd-g} we still see significant numerical slow down of the higher wave number modes.

\begin{obs}
The choice of finite difference stencil used in the coarse propagator can have a significant effect on Parareal.
It seems that centred stencils that fail to remove high wave number modes cause similar problems as non-diffusive time stepping methods, suggesting that stencils with upwind-bias are a much better choice.
\end{obs}

\subsection{Influence of phase and amplitude errors}
To investigate whether \emph{phase errors} or \emph{amplitude errors} in the coarse method trigger the instability, we construct coarse propagators where either the phase or the amplitude is exact.
Denote as $R_{\text{euler}}$ the stability function of backward Euler and as $R_{\text{exact}}$ the stability function of the exact integrator.
Then, a method with the amplification factor of backward Euler and no phase speed error can be constructed as
\begin{equation}
	\label{eq:coarse_exact_phase}
	R_{1} := \left| R_{\text{euler}}\right| e^{i \theta(R_{\text{exact}})}
\end{equation}
while a method with no amplification error and the phase speed of backward Euler can be constructed as
\begin{equation}
	\label{eq:coarse_exact_amp}
	R_{2} := \left| R_{\text{exact}} \right| e^{i \theta(R_{\text{euler}})}.
\end{equation}
These artificially constructed propagators are now used within Parareal.

\begin{figure}
	\centering
	\includegraphics[scale=1]{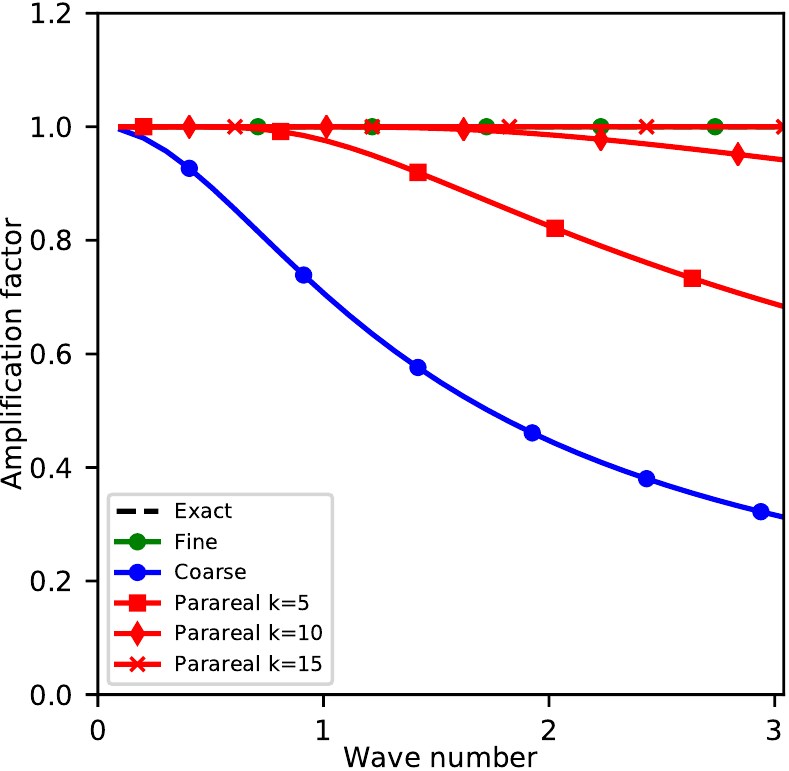}
	\includegraphics[scale=1]{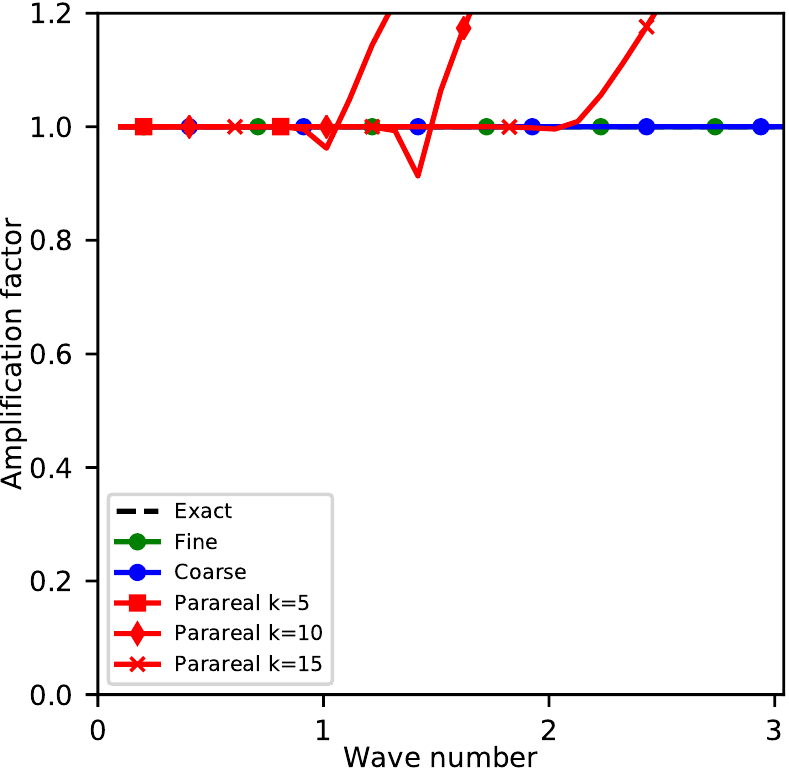}
	\caption{Amplification factor of Parareal for the advection equation for an artificially constructed coarse method with exact phase speed (upper) or exact amplification factor (lower).}
	\label{fig:phase_error_coarse}
\end{figure}

Figure~\ref{fig:phase_error_coarse} shows the resulting amplification factor when using $R_1$ (upper) or $R_2$ (lower) as coarse propagator.
For $R_1$, where there is no phase speed error in the coarse propagator, there is no instability. 
Already for $k=10$ it produces a good approximation of the exact amplification factor across the whole spectrum.
In contrast, for $R_2$ where there is no amplification error produced by $\mathcal{G}$, the instability is clearly present for $k=5$, $k=10$ and $k=15$.

\begin{figure*}[!ht]
	\centering
	\includegraphics[scale=1]{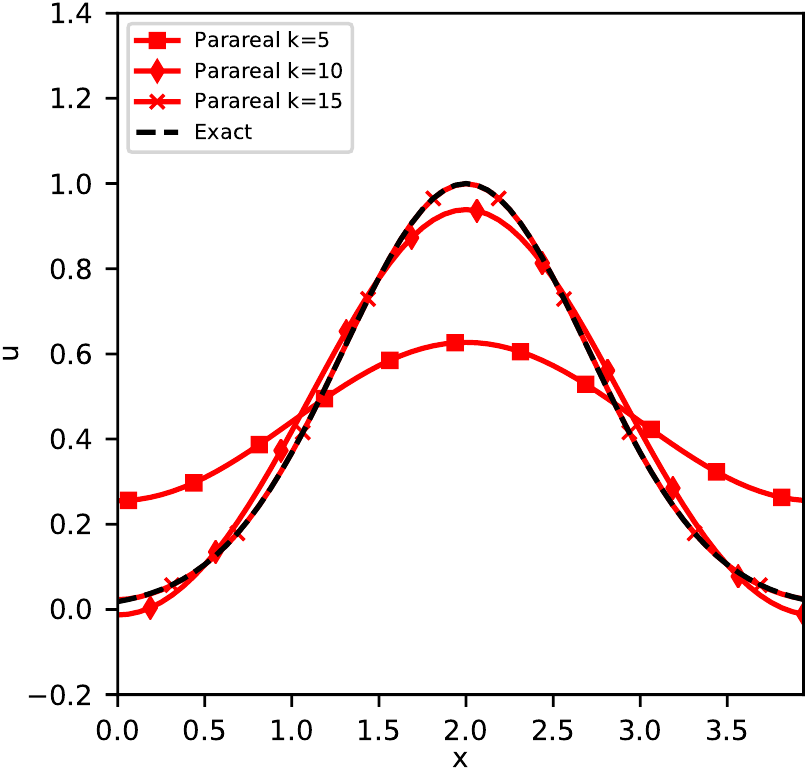}\hfill
	\includegraphics[scale=1]{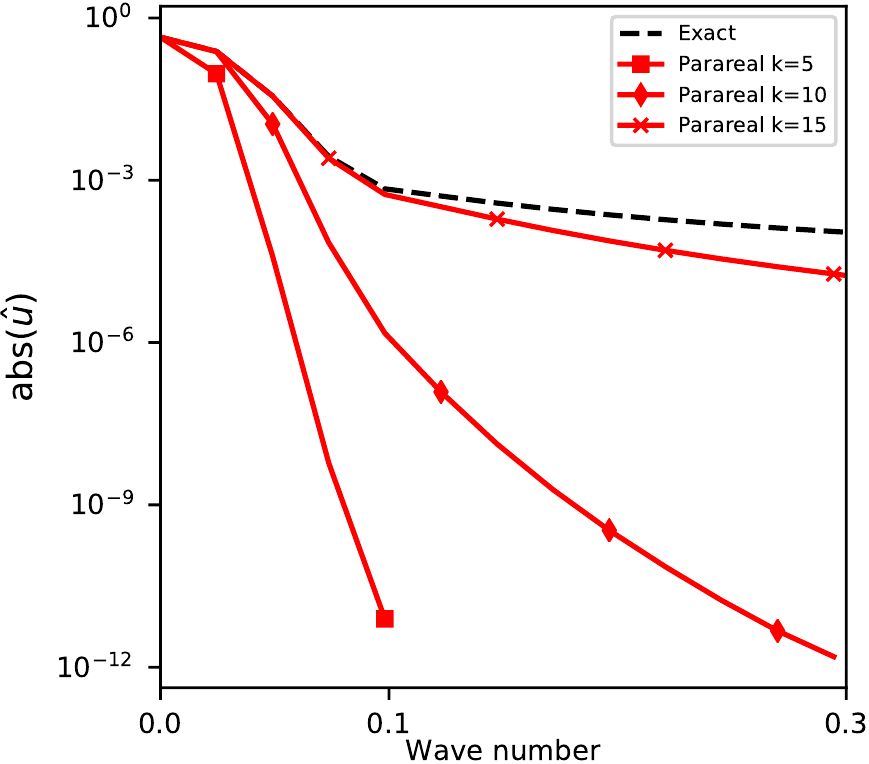}
	\caption{The same Gauss peak (left) and corresponding spectrum (right) as in Figure~\ref{fig:gauss_peak} but now computed with the $R_1$ coarse propagator with exact phase speed. }
	\label{fig:gauss_peak_exact_phase}
\end{figure*}

Figure~\ref{fig:gauss_peak_exact_phase} shows the solution for the same setup that was used for Figure~\ref{fig:gauss_peak}, except using the $R_1$ artificial coarse propagator without phase errors instead of the backward Euler.
For $k=5$ iterations, the peak is strongly damped but, because $\mathcal{G}$ has no phase errors, in the correct place.
After $k=10$ iterations, Parareal has corrected for most the numerical damping and already provides a reasonable approximation, even though the amplitude of most wave numbers in the spectrum is still severely underestimated.
However, the lack of phase errors and resulting numerical dispersion avoids the ``bulge'' and distortions that were present in Figure~\ref{fig:gauss_peak}.
Finally, for $k=15$ iterations, the solution provided by Parareal is indistinguishable from the exact solution.
Small underestimation of the amplitudes of larger wave numbers can still be seen in the spectrum, but the effect is minimal.
Note that this does not mean that Parareal will provide speedup - in a realistic scenario, where $\mathcal{F}$ is not exact but a time stepping method, too, it would depend on how many iterations are required for Parareal to be as accurate as the fine method run serially and the actual runtimes of both propagators.
All that can be said so far is that avoiding coarse propagator phase errors avoids the instability and leads to faster convergence.

\begin{figure}[!ht]
	\centering
	\includegraphics[scale=1]{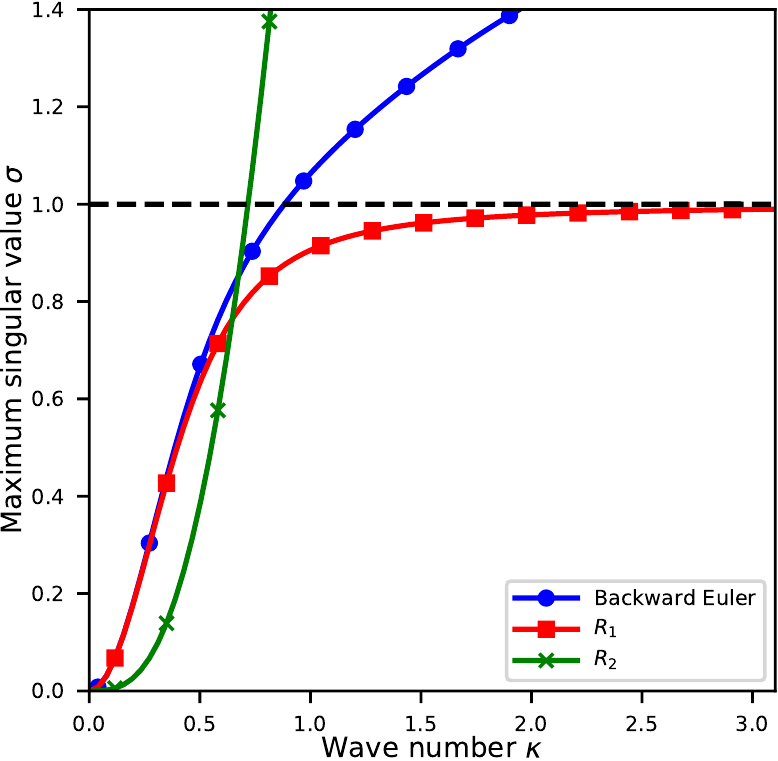}
	\caption{Maximum singular value $\sigma$ of error propagation matrix $\Tt{E}$ depending on the wave number of three choices of coarse propagator. $R_1$ has exact phase speed while $R_2$ has exact amplification factor.}
	\label{fig:sigma_vs_waveno_different_G}
\end{figure}

The effect of eliminating phase errors in the coarse method can also be illustrated by analysing the maximum singular value $\sigma$ of the error propagation matrix.
Figure~\ref{fig:sigma_vs_waveno_different_G} shows $\sigma$ depending on the wave number $\kappa$ for three different coarse propagators: the backward Euler, the artificially constructed propagator $R_1$ with no phase error and the artificially constructed propagator $R_2$ with no amplitude error.
For the backward Euler method, $\sigma$ is larger than one for significant parts of the spectrum, indicating possible non-monotonous convergence for these modes.
The situation is even worse for $R_2$, mirroring the problems with a non-diffusive coarse method like the trapezoidal rule mentioned above.
For $R_1$, however, $\sigma$ remains below one across the whole spectrum, so that Parareal will converge monotonically for every mode.
Since $\sigma$ approaches one for medium to high wave numbers, convergence there is potentially very slow, in line with the errors seen in the upper part of the spectrum of the Gauss peak.
However, in contrast to the other two cases, these wave numbers will not trigger instabilities.

In summary, these results strongly suggest that \emph{phase errors} in the coarse method are responsible for the instability, which is in line with previous findings that Parareal can quickly correct even for very strong numerical diffusion as long as a wave is placed at roughly the correct position by the coarse predictor~\cite{RuprechtKrause2012}.

\begin{obs}
The instability in Parareal seems to be caused by phase errors in the coarse propagator while amplitude errors are quickly corrected by the iteration.
\end{obs}

\subsubsection*{Relation to asymptotic Parareal}
It is interesting to point out how the $R_1$ propagator with exact phase speed is related to the asymptotic Parareal method developed by Haut et al.~\cite{HautWingate2014}.
The exact propagator for~\eqref{eq:advection_diffusion} is given by
\begin{equation}
	R_{\text{exact}} = e^{\delta(U, \kappa, \nu)} = e^{-\nu \kappa^2 t} e^{-U i \kappa}.
\end{equation}
Therefore, we have
\begin{equation}
	\left| R_{\text{exact}} \right| = e^{-\nu \kappa^2}
\end{equation}
and
\begin{equation}
	\theta(R_{\text{exact}}) = -U \kappa.
\end{equation}
Equivalent to the use of a coarse propagator $R_1$ with exact phase propagation would be solving a transformed coarse level problem instead by setting
\begin{equation}
	\tilde{u}(t) := e^{U i \kappa t} \hat{u}(t).
\end{equation}
This leads to the purely diffusive coarse level problem
\begin{equation}
	\tilde{u}_t(t) = -\nu \kappa^2 \tilde{u}(t)
\end{equation}
with restriction operator $e^{U i \kappa t}$ and interpolation $e^{- U i \kappa t}$ taking care of the propagation part.
This is precisely the strategy pursued in the nonlinear case by ``asymptotic Parareal'' where they factor out a fast term with purely imaginary eigenvalues, related to acoustic waves.
In a sense, their approach can be understood as an attempt to construct a coarse method with minimal phase speed error.
Of course, evaluation of the transformation is not trivial for more complex problems and requires a sophisticated approach~\cite{SchreiberEtAl2016}, in contrast to the here studied linear advection-diffusion equation where the transformation is simply multiplication by $e^{-U i \kappa t}$ and $e^{U i \kappa t}$.

\subsubsection*{Phase error or mismatch}
So far, we have always assumed that the fine method is exact.
This leaves the question whether the instability is triggered by phase errors in the coarse method or simply by a mismatch in phase speeds between fine and coarse level.
In order to see if the instability arises if both fine and coarse level have the same large phase error, we replace the fine propagator stability function by
\begin{equation}
	R_{3} = \left| R_{\text{fine}} \right| e^{i \theta(R_{\text{coarse}})}.
\end{equation}
Now, the fine propagator is a method with exact amplification factor but a discrete phase speed that is as inaccurate as the coarse method.
While such an integrator would not make for a very useful method in practice it is valuable for illustrative purposes.
The coarse method is again the standard implicit Euler.

\begin{figure}
	\centering
	\includegraphics[scale=1]{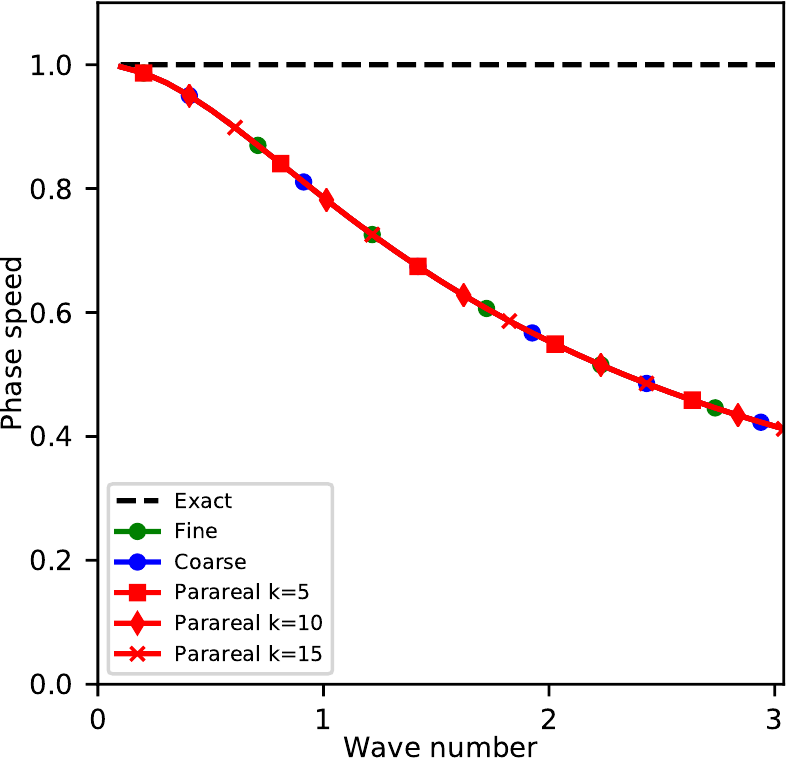}\hfill
	\includegraphics[scale=1]{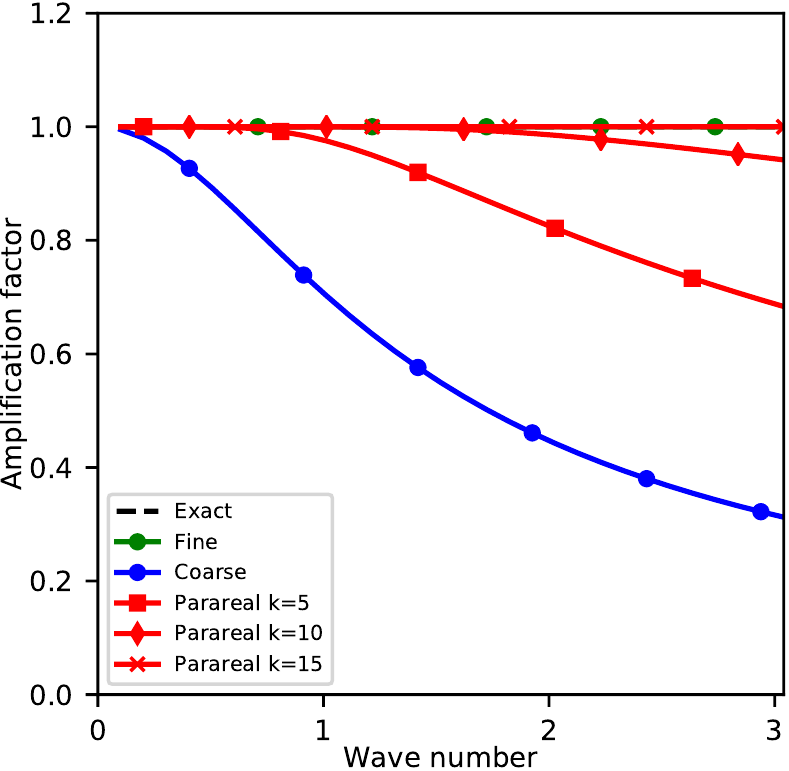}
	\caption{Phase speed (upper) and amplification factor (lower) for an artificially constructed fine propagator with the same phase error as the implicit Euler coarse propagator.}
	\label{fig:dispersion-matched-phase}
\end{figure}

Figure~\ref{fig:dispersion-matched-phase} shows the phase speed (left) and amplification factor (right) of Parareal for this configuration.
Since the fine and coarse method have the same (highly inaccurate) phase speed, Parareal matches the fine method's phase speed exactly from the first iteration and all lines coincide.
The amplification factor converges quickly to the correct value and looks almost identical to the case where a coarse propagator with exact phase speed was used, compare for Figure~\ref{fig:phase_error_coarse} (left).
No instability occurs and amplification factors are below one across the whole spectrum for all iterations.

\begin{figure*}
	\centering
	\includegraphics[scale=1]{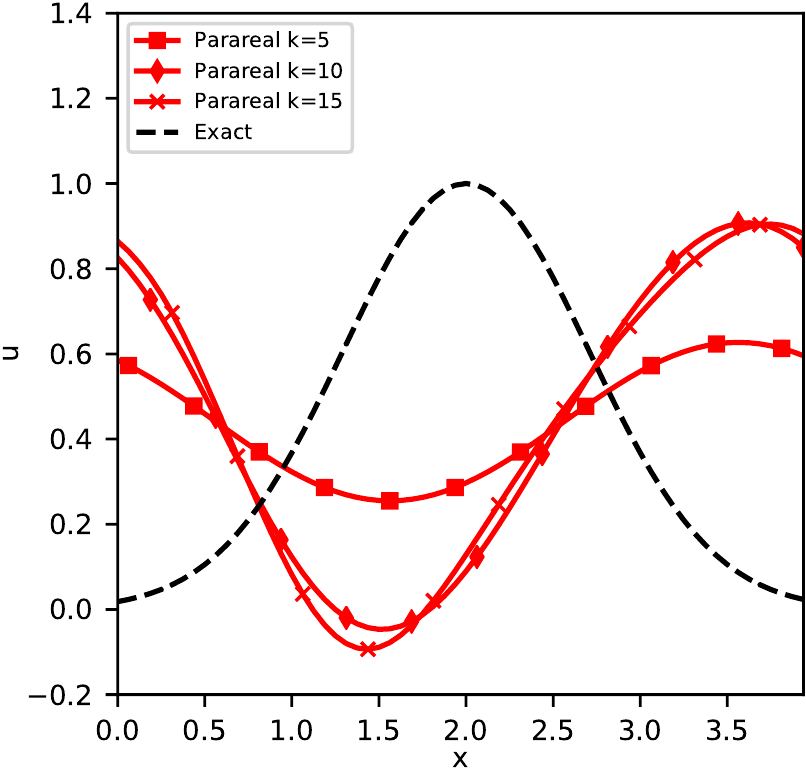}\hfill
	\includegraphics[scale=1]{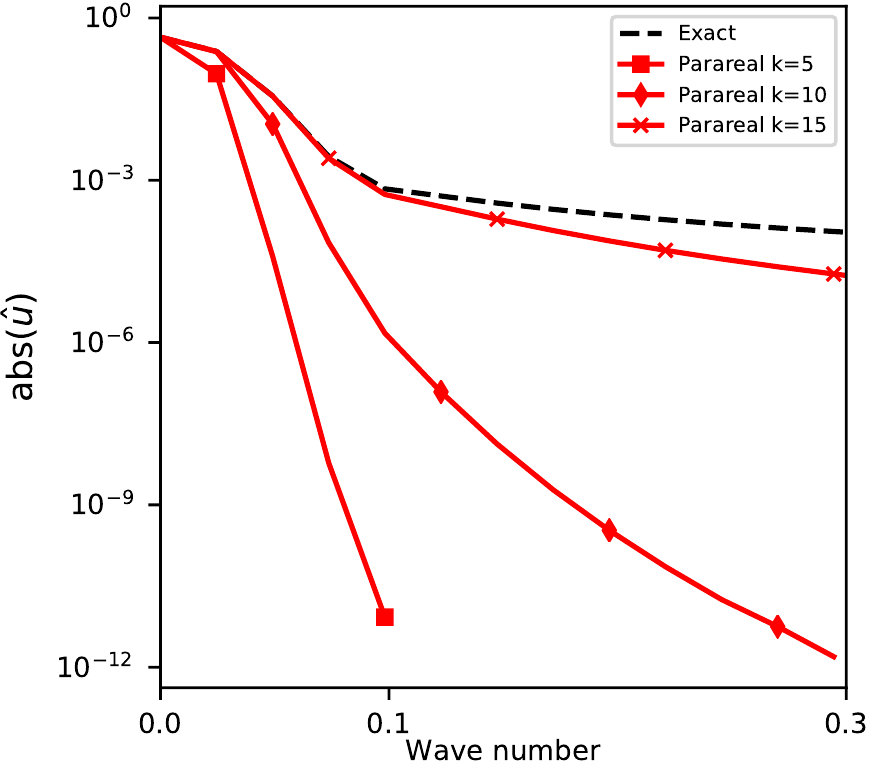}
	\caption{Gauss peak computed with the $R_3$ fine propagator.  The incorrect phase speed of the fine method puts the peak at a completely wrong position (left), but there is no instability and the spectrum (right) converges as quickly as for the exact phase speed coarse propagator in Figure~\ref{fig:gauss_peak_exact_phase}.}
	\label{fig:gauss-peak-matched-phase}
\end{figure*}

Figure~\ref{fig:gauss-peak-matched-phase} shows how Parareal converges for this configuration in physical and spectral space.
Because both fine and coarse method now have substantial phase error, the Gauss peak is at a completely wrong position.
However, for $k=10$, Parareal already approximates it reasonably well and shows no sign of instability.
Convergence looks again very similar to the results shown in Figure~\ref{fig:gauss_peak_exact_phase} except for the wrong position of the Gauss peak.
While making the fine method as inaccurate as the coarse method is clearly not a useful strategy to stabilise Parareal, this experiment nevertheless demonstrates that the instability is triggered by different discrete phase speeds in the coarse and fine method.

\begin{obs}
Analysing further the issue of phase errors shows that the instability seems to arise from mismatches between the phase speed of coarse and fine propagator.
\end{obs}

It is interesting to note that a very similar observation was made by Ernst and Gander for multi-grid methods for the Helmholtz equation.
There, the ``coarse grid correction fails because of the incorrect dispersion relation (phase error) on coarser and coarser grids [...]''~\cite{ErnstGander2013}.
They find that adjusting the wave number of the coarse level problem in relation to the mesh size leads to rapid convergence of their multi-grid solver.
Investigating if and how their approach might be applied to Parareal (which can also be considered as a multi-grid in time method~\cite{GanderVandewalle2007}) would be a very interesting direction for future research.
Furthermore, it seems possible that parallel-in-time methods with more than two levels like MGRIT~\cite{FalgoutEtAl2014_MGRIT} or PFASST~\cite{EmmettMinion2012} could yield some improvement, because they would allow for less drastic changes in resolution compared to two-level Parareal.

\subsection{Coarse time step}
\begin{figure*}[!ht]
	\centering
	\includegraphics[scale=1]{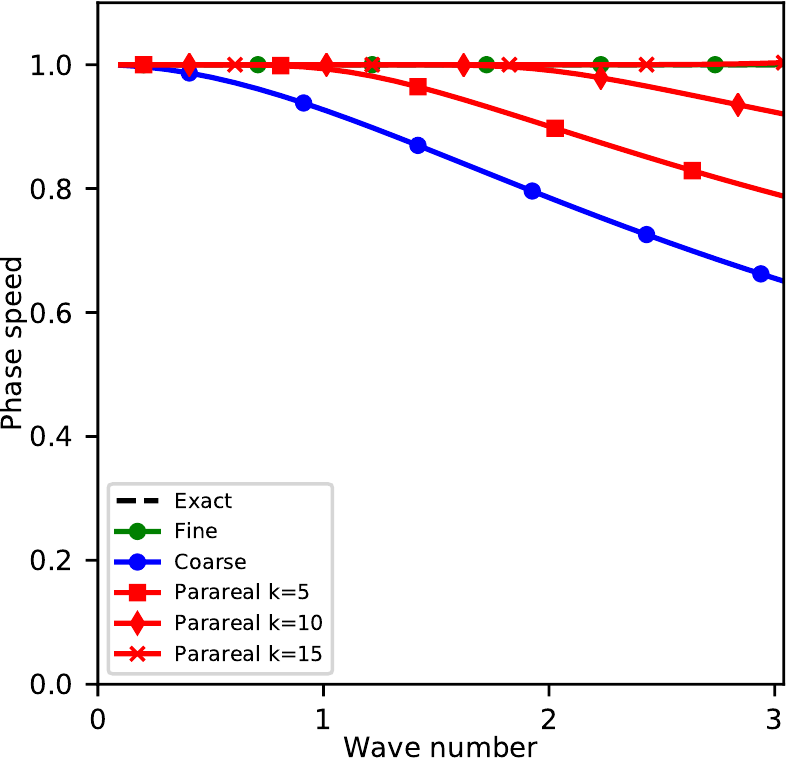}\hfill
	\includegraphics[scale=1]{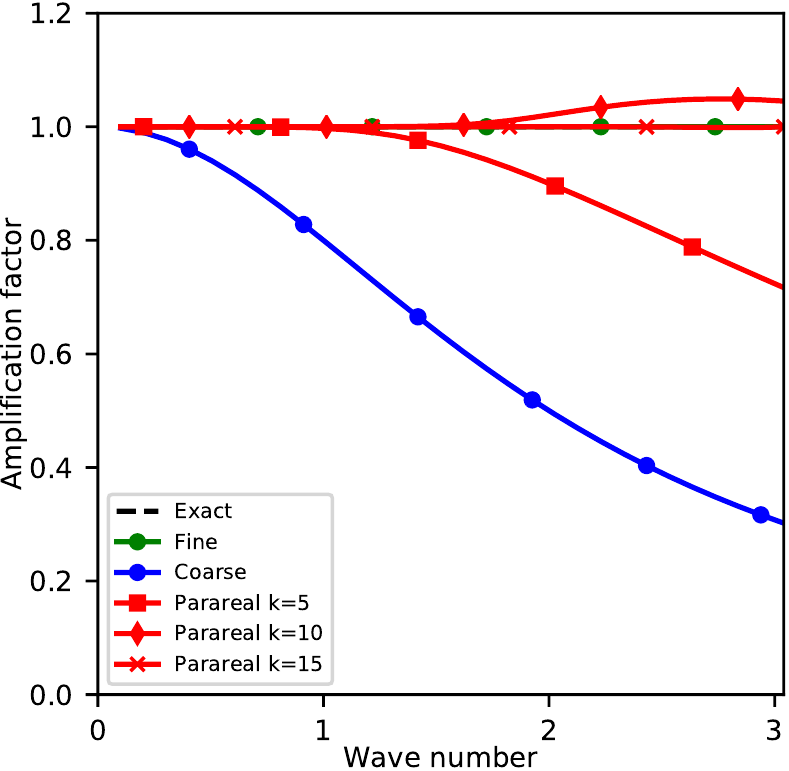}
	\caption{Phase speed (left) and amplification factor (right) for standard Parareal with the same configuration as for Figure~\ref{fig:dispersion} except using two coarse time steps per time slice.}
	\label{fig:dispersion-ncoarse-2}
\end{figure*}

Using a smaller time step for the coarse method will obviously reduce its phase error and can thus be expected to benefit Parareal convergence.
Figure~\ref{fig:dispersion-ncoarse-2} shows that this is indeed true.
It shows discrete phase speed (left) and amplification factor (right) for the same configuration as used for Figure~\ref{fig:dispersion}, except now using two coarse step per time slice instead of one.
Since the coarse propagator alone is now already significantly more accurate, Parareal with $k=5$ and $k=10$ iterations provides more accurate phase speeds and, for $k=15$, reproduces the exact value exactly,
The reduced phase errors translate into a milder instability.
For $k=10$, some wave numbers have amplification factors above one, but both the range of unstable wave numbers and the severity of the instability are much smaller than if only a single coarse time step is used.
This explains why configurations can be quite successful where both $\mathcal{F}$ and $\mathcal{G}$ use nearly identical time steps and the difference in runtime is achieved by other means, e.g. an expensive high order spatial discretisation for the fine and a cheap low order discretisation on the coarse level.

\begin{obs}
Since phase errors of the coarse method obviously depend on its time step size, reducing the coarse time step helps to reduce the range of unstable wave numbers and the severity of the instability.
\end{obs}

\subsection{Number of time slices}
\begin{figure*}[!ht]
	\centering
	\includegraphics[scale=1]{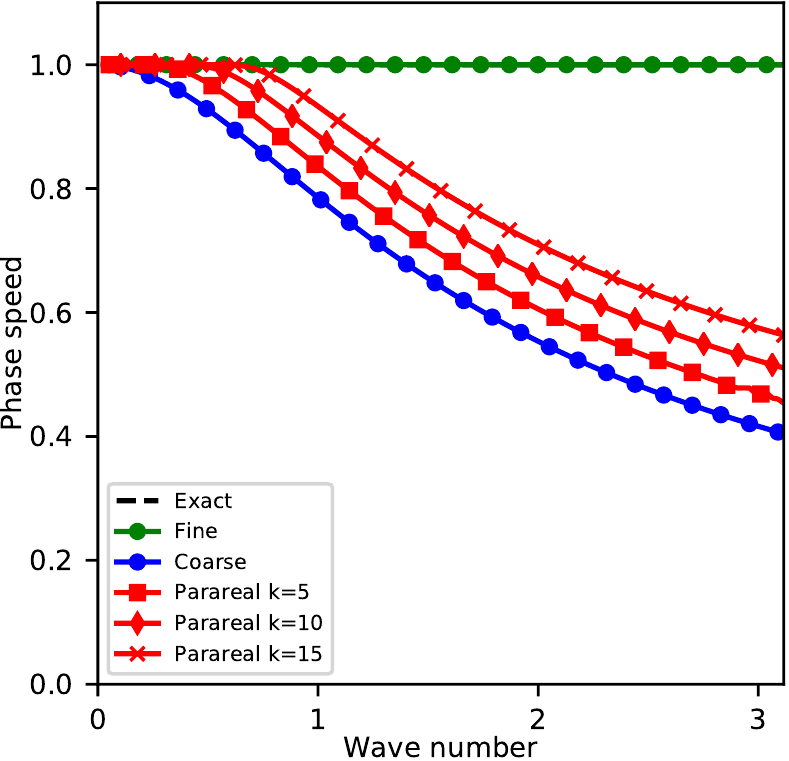}\hfill
	\includegraphics[scale=1]{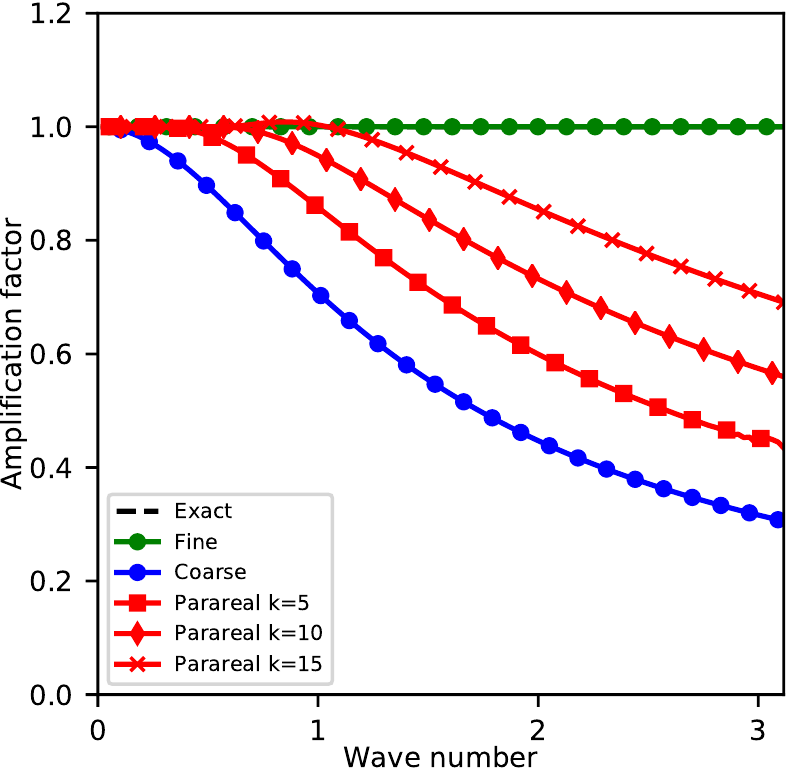}
	\caption{Phase error and amplification factor for $P=64$ in contrast to $P=16$ as in Figure~\ref{fig:dispersion}.}
	\label{fig:dispersion-P64}
\end{figure*}

All examples so far only considered $P=16$ time slices and processors.
To illustrate the effect of increasing $P$, Figure~\ref{fig:dispersion-P64} shows the discrete dispersion relation for standard Parareal for $P=64$ time slices or processors (same configuration as in Figure~\ref{fig:dispersion} except for $P$).
Even for $k=15$ iterations, Parareal reproduces the correct phase speed  (left figure) very poorly -- waves across large parts of the spectrum suffer from substantial numerical slowdown.
Also, converge is slow and there is only marginal improvement from $k=5$ to $k=15$ iterations.
Convergence is somewhat faster for the amplification factor (right figure) with more substantial improvement for $k=15$ over the coarse method.
However, there also remains significant numerical attenuation of the upper half of the spectrum.
If integrating the Gauss peak with this configuration, the result at $T=64$ after $k=5$ iterations is essentially a straight line (not shown) as almost all modes beyond $\kappa=0$ are strongly damped.
A small overshoot at around $\kappa =1$ is noticeable for $k=15$ iterations and this will worsen as $k$ increases.
In general, as $P$ increases, it takes more iterations to trigger the instability since the slow convergence requires longer to correct for the strong diffusivity of the coarse method.

\begin{figure}[!t]
	\centering
	\includegraphics[scale=1]{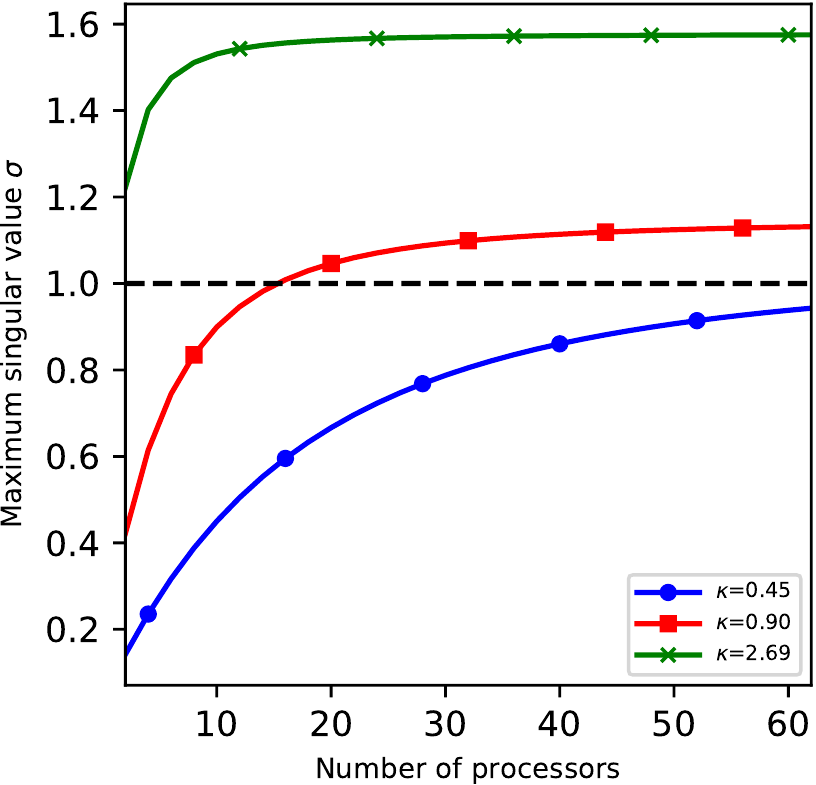}
	\caption{Maximum singular value of $\Tt{E}$ depending on the number of processors $P$ for three different wave numbers $\kappa$.}
	\label{fig:svd_vs_p}
\end{figure}

These results suggest that the high wave numbers are the slowest to converge and that convergence deteriorates as $P$ increases.
This is confirmed by Figure~\ref{fig:svd_vs_p}, showing the maximum singular value for three wave numbers plotted against $P$.
Convergence generally gets worse as $P$ increases, but note that even for $P=64$ the low wave number mode (blue circles) still converges monotonically while the high wave number mode (green crosses) might already converge non-monotonically for only $P=4$ processors.
There also seems to be a limit for $\sigma$ as $P$ increases, with higher wave numbers levelling off at higher values of $\sigma$.

Therefore, Parareal could provide some speedup for linear hyperbolic problems if the solution consists mainly of very low wave number modes and/or numerical diffusion in the fine propagator is sufficiently strong to remove higher wave number modes.
This also explains why divergence damping in the fine propagator can accelerate convergence of Parareal~\cite{RuprechtKrause2012}, as it removes exactly the high wave number modes that converge the slowest.

\begin{obs}
While convergence becomes slower as the number of processors $P$ increases, low wave numbers converge monotonically even for large numbers of $P$ but high wave numbers might not do so already for $P=4$.
\end{obs}

\subsection{Wave number}
\begin{figure}[!t]
	\centering
	\includegraphics[scale=1]{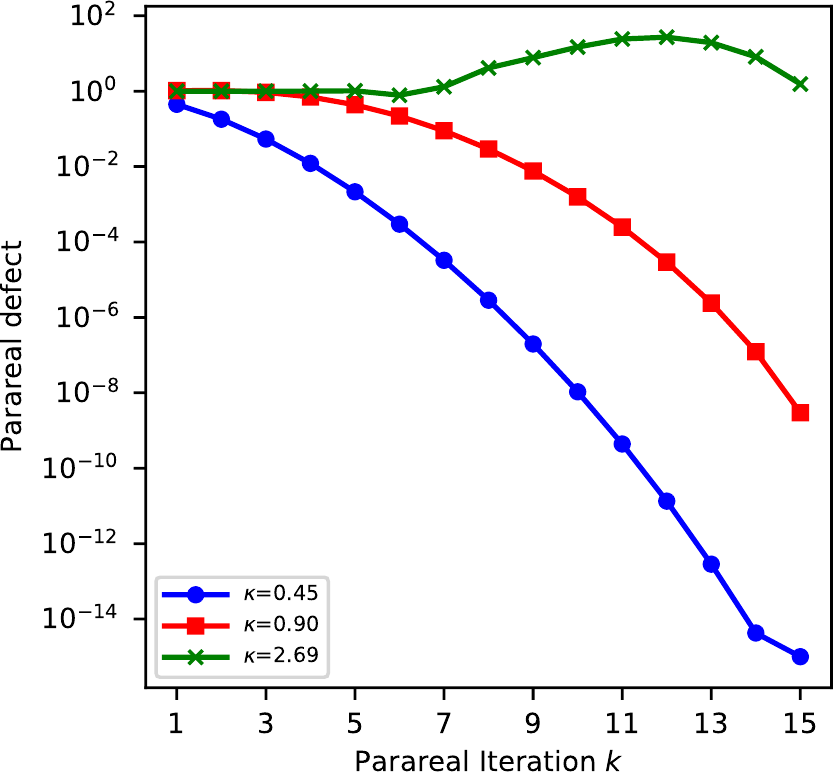}\hfill
	\includegraphics[scale=1]{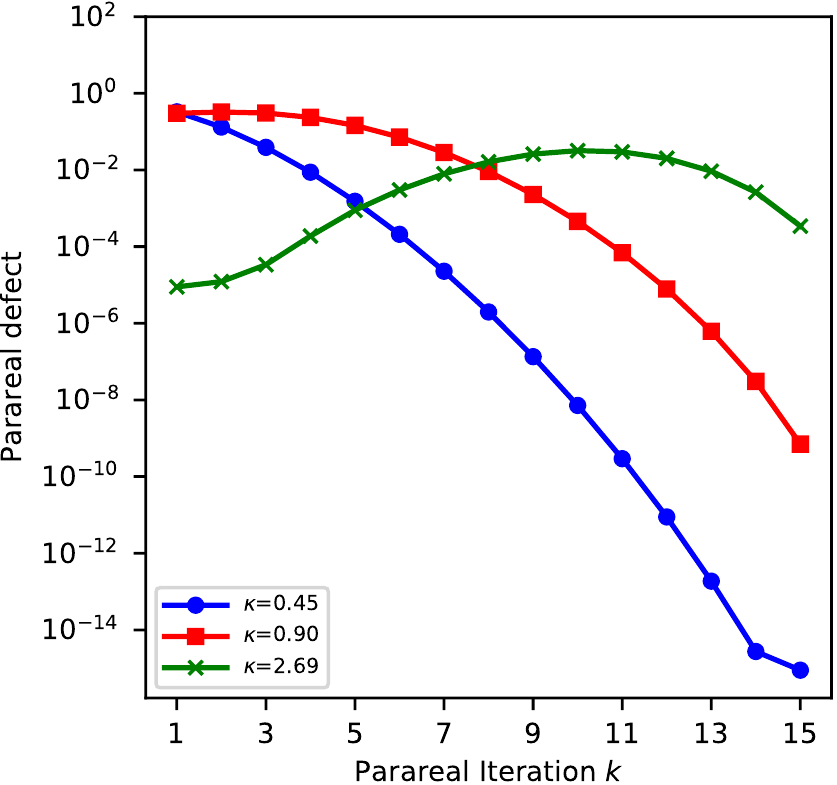}
	\caption{Parareal defect versus number of iterations for $\nu=0.0$ (upper) and $\nu=0.1$ (lower).}
	\label{fig:conv_waveno}
\end{figure}

The analysis above showed that higher wave numbers converge slower and are more susceptible to instabilities.
This is confirmed in Figure~\ref{fig:conv_waveno} showing the difference between the Parareal and fine integrator stability function
\begin{equation}
	d(k) := \left| R_{\text{parareal}}(k) - R_{\text{fine}} \right|.
\end{equation}
The smallest wave number, $\kappa = 0.45$, converges quickly in the hyperbolic (upper) and diffusive case (lower).
For $\kappa=0.9$, both cases show some initial stalling before the mode converges.
Finally, $\kappa = 2.69$ shows first a significant increase in the defect before convergence sets in as $k$ comes close to $P=16$.
Interestingly, this is the case for both $\nu=0$ and $\nu=0.1$. 
While the ``bulge'' is even more pronounced in the diffusive case, since modes decay proportional to $e^{-\nu \kappa^2 t}$ in amplitude, the absolute values of the defect are orders of magnitude smaller.
Therefore, at least until $k=7$ iterations, it is no longer the high wave number mode $\kappa = 2.69$ that will restrict performance, but rather the lower wave number $\kappa = 0.9$.
Then, the instability for the high wave number kicks in and for $k \geq 8$ wave number $\kappa = 2.69$ is again causing the largest defect.
As $\nu$ increases, however, the defects for $\kappa = 2.69$ will reduce further, the cross-over point will move to later iterations and eventually lower wave numbers will determine convergence for all iterations.
In a sense, in line with the analysis above, Parareal propagates high wave number modes very wrongly in both cases, but since high wave number modes are quickly attenuated if $\nu > 0$, it does not matter very much in the diffusive case.

\begin{figure}[!t]
	\centering
	\includegraphics[scale=1]{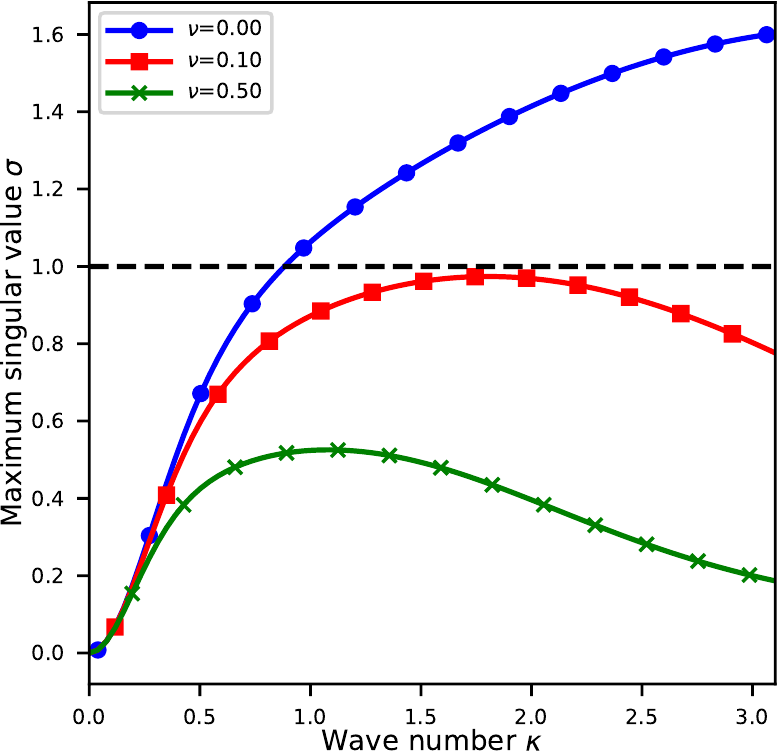}
	\caption{Maximum singular value of $\Tt{E}$ depending on the wave number $\kappa$ for three values of diffusivity $\nu$.}
	\label{fig:svd_vs_waveno}
\end{figure}

Figure~\ref{fig:svd_vs_waveno} confirms this for a wider range of wave numbers $\kappa$.
It shows the maximum singular value $\sigma$ of the error propagation matrix $\Tt{E}$ over the whole spectrum for three different values of $\nu$.
For $\nu = 0$ (hyperbolic case), $\sigma$ increases monotonically with $\kappa$ and the highest wave number converges the slowest.
After around $\kappa \geq 0.8$, the singular values are larger than one and convergence becomes potentially non-monotone.
For $\nu = 0.1$, $\sigma$ increases until around $\kappa = 1.8$ and then decreases again.
Therefore, the slowest converging mode is no longer at the end but in the middle of the spectrum.
Also, we now have $\sigma < 1$ for all $\kappa$ so that all modes will converge, even though some potentially very slowly.
Increasing diffusion further to $\nu = 0.5$ greatly improves convergence for all modes, the largest $\sigma$ across the whole spectrum is now below $0.5$.
The worst converging mode has also moved ``further down'' the spectrum and is now at around $\kappa = 1.0$.
This shows how the strong natural damping of high wave numbers from diffusion counteracts Parareal's tendency to amplify them and thus stabilises it.

\begin{obs}
Since diffusion naturally damps higher wave numbers, it removes the issue of slow or no convergence at the end of the spectrum.
Therefore, as the diffusivity parameter $\nu$ increases, the wave number that converges the slowest and determines performance of Parareal becomes smaller.
\end{obs}

\section{Conclusions}
Efficient parallel-in-time integration of hyperbolic and ad\-vection-dominated problems has been shown to be problematic.
This prevents application of a promising new parallelisation strategy to many problems in computational fluid dynamics, despite the urgent need for better utilisation of massively parallel computers.
For the Parareal parallel-in-time method, mathematical theory has shown that the algorithm is either unstable or inefficient when applied to hyperbolic equations, but so far no detailed analysis exists of how exactly the instability emerges.

The paper presents the first detailed analysis of how Parareal propagates waves and the ways in which the instability is triggered.
It uses a formulation of Parareal as a preconditioned fixed point iteration for linear problems to derive its stability function.
From there, a discrete dispersion relation is obtained that allows to study the phase speed and amplitude errors from Parareal when computing wave-like solutions.
To deal with issues arising from increasing the time interval together with the number of processors, a simple procedure is introduced to normalise the stability function to the unit interval.

Analysis of the discrete dispersion relation and the maximum singular value of the error propagation matrix then allows to make a range of observations, illustrating where the issues of Parareal for wave problems originate.
A key finding is that the source of the instability are different discrete phase speeds on the coarse and fine level, which cause instability of higher wave number modes.
Interestingly, the overestimation of high wave number amplitudes is present in diffusive problems, too, but since there these amplitudes are naturally strongly damped, this does not trigger instabilities.
Further analysis addresses the role of the number of processors, the coarse time step size and comments on possible connections to asymptotic Parareal and multi-grid methods for the Helmholtz equation.

The analysis presented here will be useful to interpret and understand performance of Parareal for more complex problems in computational fluid dynamics.
A natural line of future research would be to attempt to develop a new, more stable, parallel-in-time method for hyperbolic problems based on the provided observations.
For example, the update in Parareal proceeds component wise.
That means that if the coarse propagator moves a wave at the wrong speed, the update will not know that a simple shift of entries could provide a good correction.
Attempting to somehow modify the Parareal update to take into account this type of information seems promising, even though probably challenging to do in 3D.
Extending the analysis presented here to systems with multiple waves, e.g. the shallow water equations, or to nonlinear problem where wave numbers interact would be another interesting line of inquiry.
Furthermore, the framework used here to analyse Parareal is straightforward to adopt for other parallel-in-time integration methods as long as a matrix representation for them is available.

\paragraph{Acknowledgments.} I would like to thank Martin Gander, Martin Schreiber and Beth Wingate for the very interesting discussions at the 5th Workshop on Parallel-in-time integration at the Banff International Research Station (BIRS) in November 2016, which led to the comments about asymptotic Parareal and multi-grid for 1D Helmholtz equation in the paper.

Parts of the pseudo-spectral code used to illustrate the effect of numerical dispersion come from David Ketcheson's short course PseudoSpectralPython~\cite{ps_py2015}.

The pyParareal code written for this paper relies heavily on the open Python packages \texttt{NumPy}~\cite{numpy}, \texttt{SciPy}~\cite{scipy} and \texttt{Mat\-plotlib}~\cite{matplotlib}. 

\bibliographystyle{spmpsci}      
\bibliography{pint,refs,py}   

\end{document}